\documentclass[leqno]{amsart}
\usepackage{amssymb}
\usepackage{amsfonts}
\usepackage{amsfonts}
\usepackage{amsfonts}
\usepackage{amsfonts}
\usepackage{mathrsfs}
\usepackage{amssymb}
\usepackage{txfonts}
\usepackage{amsfonts}
\usepackage{amssymb,amsmath}

\usepackage{graphicx}
\usepackage{subfigure}

\usepackage{graphicx,epic,eepic}
\usepackage{ulem}

\allowdisplaybreaks

\begin{document}

\title[]
{Bubble solution for the critical Hartree equation in pierced domain}

\author[M.G. Ghimenti]{Marco Gipo Ghimenti}
\address[Marco Ghimenti]{Dipartimento di Matematica
Universit\`a di Pisa
Largo Bruno Pontecorvo 5, I - 56127 Pisa, Italy}
\email{marco.ghimenti@unipi.it }

\author[X. Huang]{Xiaomeng Huang}
\address[Xiaomeng Huang]{School of Mathematics and Statistics, Southwest University,
Chongqing 400715, People's Republic of China.}\email{hhuangxiaomeng@126.com}

\author[A. Pistoia]{Angela Pistoia}
\address[Angela Pistoia] {Dipartimento SBAI, Universit\`{a} di Roma ``La Sapienza", via Antonio Scarpa 16, 00161 Roma, Italy}
\email{angela.pistoia@uniroma1.it}

 \maketitle

\maketitle
\numberwithin{equation}{section}
\newtheorem{theorem}{Theorem}[section]
\newtheorem{lemma}[theorem]{Lemma}
\newtheorem{definition}[theorem]{Definition}
\newtheorem{proposition}[theorem]{Proposition}
\newtheorem{remark}[theorem]{Remark}
\allowdisplaybreaks

\maketitle

\noindent {\bf Abstract}: In this article, we establish the existence of solutions to the following critical Hartree equation
\begin{align*}
\begin{cases}
-\Delta u=\left(\int_{\Omega_\varepsilon}\frac{u^{2_{\mu}^*}}{|x-y|^{\mu}}dy\right)u^{2_{\mu}^*-1}, &\text{ in } \Omega_\varepsilon, \\
u=0, &\text{ on } \partial\Omega_\varepsilon,
\end{cases}
\end{align*}
where $2_{\mu}^*=\frac{2N-\mu}{N-2}$ is the upper critical exponent in the sense of the Hardy-Littlewood-Sobolev inequality, $N\geq 5$, $0<\mu<4$ with $\mu$ sufficiently close to $0$, $\Omega_\varepsilon:=\Omega\backslash B(0,\varepsilon)$ and $\Omega$ is a bounded smooth domain in $\mathbb{R}^N$, which contains the origin, and $\varepsilon$ is a positive parameter. As $\varepsilon$ goes to zero, we construct bubble solution which blows up at the origin.

\vspace{3mm} \noindent {\bf Keywords}: Critical Hartree equation; Bubble solution; Reduction arguments.

\vspace{3mm}

\maketitle

\section{Introduction and statement of main result}\label{sec preliminaries}

In this work, we study the following critical Hartree equation
\begin{align}\label{crihartree}
\begin{cases}
-\Delta u=\left(\int_{\Omega_\varepsilon}\frac{u^{2_{\mu}^*}}{|x-y|^{\mu}}dy\right)u^{2_{\mu}^*-1}, &\text{ in } \Omega_\varepsilon, \\
u=0, &\text{ on } \partial\Omega_\varepsilon,
\end{cases}
\end{align}
where $N\geq 5$, $0<\mu<4$ with $\mu$ sufficiently close to $0$, $\Omega_\varepsilon:=\Omega\backslash B(0,\varepsilon)$ and $\Omega$ is a bounded smooth domain in $\mathbb{R}^N$, $0\in \Omega$, $2_{\mu}^*=\frac{2N-\mu}{N-2}$, is the upper critical exponent in the sense of the Hardy-Littlewood-Sobolev inequality, which we recall hereafter.

\begin{proposition}\cite{lieb,liebl}\label{hls}
Let $\theta$, $r>1$ and $0<\mu<N$ with $\frac{1}{\theta}+\frac{\mu}{N}+\frac{1}{r}=2$. Let $f\in L^{\theta}(\mathbb{R}^N)$ and $g\in L^{r}(\mathbb{R}^N)$, there exists a sharp constant $C(\theta,r,\mu,N)$ independent of $f$, $g$, such that
\begin{align}\label{hlsinequality}
\int_{\mathbb{R}^N}\int_{\mathbb{R}^N}\frac{f(x)g(y)}{|x-y|^{\mu}}dxdy
\leq C(\theta,r,\mu,N)\|f\|_{\theta}\|g\|_r.
\end{align}
If $\theta=r=\frac{2N}{2N-\theta}$, then
\begin{align*}
C(\theta,r,\mu,N)=C_{N,\mu}
=\pi^{\frac{\mu}{2}}\frac{\Gamma(\frac{N-\mu}{2})}{\Gamma(N-\frac{\mu}{2})}
\Bigg(\frac{\Gamma(N)}{\Gamma(\frac{N}{2})}\Bigg)^{\frac{N-\mu}{N}}.
\end{align*}
There is equality in (\ref{hlsinequality}) if and only if $f\equiv(const.)g$ and
\begin{align*}
g(x)=\tilde{A}(\tilde{\gamma}^2+|x-\tilde{a}|^2)^{-\frac{2N-\mu}{2}}
\end{align*}
for some $\tilde{A}\in \mathbb{C}$, $0\neq \tilde{\gamma}\in \mathbb{R}$ and $\tilde{a}\in \mathbb{R}^N$.
\end{proposition}
It follows from the Hardy-Littlewood-Sobolev inequality that
\begin{align*}
\left(\int_{\mathbb{R}^N}\int_{\mathbb{R}^N}\frac{|u(x)|^{\frac{2N-\mu}{N-2}}|u(y)|^{\frac{2N-\mu}{N-2}}}{|x-y|^{\mu}}dxdy
\right)^{\frac{N-2}{2N-\mu}}
\leq C_{N,\mu}^{\frac{N-2}{2N-\mu}}\|u\|_{2^*}^2.
\end{align*}

Considering the following critical Hartree equation
\begin{align}\label{crihartreern}
-\Delta u=
\left(\int_{\mathbb{R}^N}\frac{u^{2_{\mu}^*}}{|x-y|^{\mu}}dy\right)u^{2_{\mu}^*-1}  \ &\text{ in } \mathbb{R}^N,
\end{align}
which is related to the nonlinear Choquard equation
\begin{align}\label{nonchoquard}
-\Delta u+V(x)u=\left(\int_{\mathbb{R}^N}\frac{|u|^{p}}{|x-y|^{\mu}}dy\right)|u|^{p-2}u  \ &\text{ in } \mathbb{R}^N.
\end{align}
This equation has appeared in the context of various physical models, such as
Fr\"{o}hlich and Pekar's model of the polaron, see \cite{frohlich,frohlich1,pekar}, or the modelling of a one-component plasma introduced by Ph. Choquard \cite{Choquard}.
In the past few years, there is much work related to (\ref{nonchoquard}). If $p=\frac{2N-\mu}{N}$, Moroz and Van Schaftingen \cite{morozv1} proved the existence of one nontrivial solution to (\ref{nonchoquard}) if $V(x)$ satisfies some condition. For the case $p\in (\frac{2N-\mu}{N},\frac{2N-\mu}{N-2})$, there are many papers studying the existence, symmetry and qualitative
properties of solution to (\ref{nonchoquard}) by variational methods, see for example \cite{alvesny, maz, miaowx, morozpt, morozv}. Wei and Winter \cite{weiw} constructed families of solutions for the singularly perturbed subcritical Choquard equation by using a Liapunov-Schmidt type reduction. Subsequently,
Luo, Peng and Wang \cite{luopw} applied a local Pohozaev type of identity, blow-up analysis and
the maximum principle to prove the uniqueness of the solutions obtained by Wei and Winter \cite{weiw}.
By applying a finite dimensional reduction argument and
developing novel local Poho\v{z}aev identities, Gao, Moroz, Yang and Zhao \cite{gaomyz} proved that the problem (\ref{nonchoquard}) admits infinitely many solutions with arbitrary large energies when $N=6$ and $p=\frac{2N-\mu}{N-2}$.

Concerning (\ref{crihartreern}), by using the moving plane methods in integral form, Lei \cite{lei}, Du and Yang \cite{duy},
Guo, Hu, Peng, Shuai \cite{guohps} classified independently the positive
solutions of the critical Hartree equation (\ref{crihartreern}) and
proved that each positive solution of (\ref{crihartreern}) is radially symmetric, monotone decreasing about some point and has the form
\begin{align}\label{utx}
U_{t,\xi}(x):=c_\mu\frac{t^{\frac{N-2}{2}}}{\left(t^2+|x-\xi|^2\right)^{\frac{N-2}{2}}},\ t>0,\ \xi\in\mathbb R^N,
\end{align}
where $c_\mu>0$ depends only on $\mu$ and $N$. For the nondegeneracy result for the critical Hartree equation (\ref{crihartreern}), Du and Yang \cite{duy} considered the nondegeneracy of the unique solutions for the case of $\mu$ close to $N$, while Giacomoni, Wei and Yang \cite{giacononiwy} studied the case of $\mu$ close to $0$. Recently, Gao, Moroz, Yang and Zhao \cite{gaomyz} proved that a nondegeneracy result for the critical Hartree equation in dimension $6$. After summarizing the above nondegeneracy results, Yang and Zhao \cite{yangz} proved via a reduction method that, as $\varepsilon\rightarrow 0$, the solution to the following Hartree type Brezis-Nirenberg problem
\begin{align}\label{crihartreebn}
\begin{cases}
-\Delta u=\left(\int_{\Omega_\varepsilon}\frac{u^{2_{\mu}^*}}{|x-y|^{\mu}}dy\right)u^{2_{\mu}^*-1}+\varepsilon u, &\text{ in } \Omega, \\
u=0, &\text{ on } \partial\Omega,
\end{cases}
\end{align}
blows up exactly at a critical point of the Robin function that cannot be on the boundary of $\Omega$.
Recently, Yang, Ye and Zhao \cite{yangyz} obtained a kind of converse result obtained in \cite{yangz}, that is, they proved that problem (\ref{crihartreebn}) has a family of solutions $u_\varepsilon$ concentrating around the critical point of Robin function.
By local Pohozaev identities and blow-up analysis, Squassina, Yang and Zhao \cite{squassinayz} not only located the blow-up point of single bubbling solutions for (\ref{crihartreebn}), but also obtained the local uniqueness of the blow-up solutions for (\ref{crihartreebn}).
Using the techniques in \cite{brezisn}, Gao and Yang \cite{gaoy} extended the well-known results established by Brezis and Nirenberg in \cite{brezisn} for local elliptic equation to the nonlocal Choquard equation, more precisely, they proved that
problem (\ref{crihartreebn}) has a solution for $\varepsilon>0$ when $N\geq 4$; if $N=3$, then there exists $\varepsilon^*$ such that (\ref{crihartreebn}) has a nontrivial solution for $\varepsilon> \varepsilon^*$, where $\varepsilon$ is not an eigenvalue of $-\Delta$ with homogeneous Dirichlet boundary data; if $N\geq 3$ and $\varepsilon\leq 0$, (\ref{crihartreebn}) admits no solutions when $\Omega$ is star-shaped.
Ghimenti and Pagliardini \cite{ghimentip} considered a slightly subcritical Choquard problem on a bounded domain and got that the number of positive solutions depends on the topology of the domain.

As far as we know, the problem of a critical Choquard equation with the pierced domain has not been studied in literature yet. So,
inspired by the work \cite{yangyz,yangz}, we are interested in constructing single peak to the critical Choquard equation when the domain has a circular hole which shrinks to a point. In particular, we will prove that problem (\ref{crihartree}) has a solution $u_\varepsilon$ concentrating around the center of the hole.

The main result is as follows.

\begin{theorem}\label{thm1.1} Assume $N\geq 5$, $0<\mu<4$ with $\mu$ sufficiently close to $0.$ Then there exists $\varepsilon_0>0$ such that for any $\varepsilon\in (0, \varepsilon_0)$ there exists a solution $u_\varepsilon$ to problem (\ref{crihartree}) such that
\begin{align*}
u_\varepsilon(x)=\left(\frac{M\varepsilon^{-\frac{1}{2}}}
{1+M^2\varepsilon^{-1}|x|^2}\right)^{\frac{N-2}{2}}
+\phi_\varepsilon(x),
\end{align*}
where $M$ is a positive constant depending only on $N$ and $\|\phi_\varepsilon\|\rightarrow 0$ as $\varepsilon\rightarrow 0$.
\end{theorem}

The proof of Theorem \ref{thm1.1} follows along the general lines of that one devised for the construction of a single spike, see \cite{yangyz,yangz} and the reference therein. However the existence of the nonlocal nonlinearity with critical exponent will lead to some difficulties, and we need more fine estimate to overcome this difficulties.
We point out that the assumption on the smallness of $\mu$ is only required to ensure the non-degeneracy of the bubble \eqref{utx} (see \cite{yangyz}). Moreover, the assumption $0<\mu<4$ seems to be necessary to have a good first order approximation of the solution we are going to build (see Lemma \ref{IVphic1}). It would be interesting to check how the  ansatz could be improved when $\mu>4.$ \\

The paper is organized as follows. In Section~\ref{sec preliminaries}, we give some definitions and preliminaries.
Section \ref{Ansatz and sketch of the proof} is devoted to prove Theorem \ref{thm1.1}, we obtain the existence of solution to problem (\ref{crihartree}). The finite dimensional reduction is performed in Section \ref{The Linear Problem}, while Section \ref{energy expansion} contains the asymptotic expansion of the reduced energy. The rest of the paper is devoted to useful estimate used to reach the results in Theorem \ref{thm1.1}.


\section{Ansatz and sketch of the proof}\label{Ansatz and sketch of the proof}

In order to state our main result it is necessary to introduce the bubble, which is the key ingredient of our proof. A bubble is the function $U_{\lambda_{\varepsilon},\xi_{\varepsilon}}$ defined by
\begin{align*}
U_{\lambda_{\varepsilon},\xi_{\varepsilon}}(x)=\frac{
\lambda_{\varepsilon}^{\frac{N-2}{2}}}{\left(1+\lambda_{\varepsilon}^2|x-\xi_{\varepsilon}|^2\right)^{\frac{N-2}{2}}},\ x, \xi_{\varepsilon}\in \mathbb{R}^N, \lambda_{\varepsilon}\in \mathbb{R}^+.
\end{align*}
It is well known (see \cite{aubin,talenti}) that they are all the positive solutions of the limit problem
\begin{align*}
\begin{cases}
-\Delta u=N(N-2)u^{2^*-1} & \text{ in } \mathbb{R}^N\\
u>0   & \text{ in } \mathbb{R}^N\\
u\in H^1(\mathbb{R}^N),
\end{cases}
\end{align*}
that achieves the best Sobolev constant $S$. Denote by
\begin{align}\label{AHL}
\mathcal{A}_{H,L}=\left(N(N-2)\right)^{\frac{N-\mu+2}{2}}\frac{S^{\frac{\mu-N}{2}}}{C(N,\mu)},
\end{align}
where $C(N,\mu)$ is the sharp constant of Hardy-Littlewood-Sobolev inequality given as in (\ref{hlsinequality}).
Then $U_{\lambda_{\varepsilon},\xi_{\varepsilon}}$ satisfies the equation
\begin{align}\label{Ulimiteq}
-\Delta U_{\lambda_{\varepsilon},\xi_{\varepsilon}}=\mathcal{A}_{H,L}\left(\frac{1}{|x|^{\mu}}\ast U^{2^*_{\mu}}_{\lambda_{\varepsilon},\xi_{\varepsilon}}\right)U^{2^*_{\mu}-1}_{\lambda_{\varepsilon},\xi_{\varepsilon}} \ \text{ in } \mathbb{R}^N.
\end{align}
Moreover, $U_{\lambda_{\varepsilon},\xi_{\varepsilon}}$
are the only solutions of (\ref{Ulimiteq}), and the kernel of the linearized operators
\begin{align*}
L\varphi:=-\Delta \varphi
-2_{\mu}^*\left(\frac{1}{|x|^{\mu}}\ast(U_{\lambda_{\varepsilon},\xi_{\varepsilon}}^{2_{\mu}^*-1})\varphi\right)
U_{\lambda_{\varepsilon},\xi_{\varepsilon}}^{2_{\mu}^*-1}
-(2_{\mu}^*-1)\left(\frac{1}{|x|^{\mu}}\ast U_{\lambda_{\varepsilon},\xi_{\varepsilon}}^{2_{\mu}^*}\right)U_{\lambda_{\varepsilon},\xi_{\varepsilon}}^{2_{\mu}^*-2}\varphi
\end{align*}
around these solutions is generated by
\begin{align*}
Z^0_{\lambda_{\varepsilon},\xi_{\varepsilon}}=\frac{\partial U_{\lambda_{\varepsilon},\xi_{\varepsilon}}}{\partial \lambda_\varepsilon}
=\frac{N-2}{2}\lambda_\varepsilon^{\frac{N-4}{2}}\frac{1-\lambda_\varepsilon^2|x-\xi_\varepsilon|^2}
{(1+\lambda_\varepsilon^2|x-\xi_\varepsilon|^2)^{\frac{N}{2}}},
\end{align*}
\begin{align*}
Z^j_{\lambda_{\varepsilon},\xi_{\varepsilon}}=\frac{\partial U_{\lambda_{\varepsilon},\xi_{\varepsilon}}}{\partial \xi^j_\varepsilon}
=(N-2)\lambda_\varepsilon^{\frac{N+2}{2}}\frac{(x-\xi_\varepsilon)_j}{(1+\lambda_\varepsilon^2|x-\xi_\varepsilon|^2)^{\frac{N}{2}}}, \ j=1,\cdots,N,
\end{align*}
(see \cite{yangyz}). Nondegeneracy is a fundamental requirement to perform Liapunov-Schmidt reduction. In fact, this is
needed to prove Lemma \ref{Linver}, which is the first step of the finite dimensional reduction.

To find a good approximation of the solution we are looking for, we need to project the bubble onto the domain $\Omega_\varepsilon$ with Dirichlet boundary conditions. For any function $U_{\lambda_{\varepsilon},\xi_{\varepsilon}}\in D^{1,2}(\mathbb{R}^N)$ we denote by $PU_{\lambda_{\varepsilon},\xi_{\varepsilon}}$ its protection on $H_0^1(\Omega_\varepsilon)$, i.e. the unique solution of the problem
\begin{align}\label{PUeq}
\begin{cases}
-\Delta \omega=\mathcal{A}_{H,L}\left(\int_{\Omega_\varepsilon}\frac{U^{2^*_{\mu}}_{\lambda_{\varepsilon},\xi_{\varepsilon}}(y)}{|x-y|^{\mu}}dy \right)U^{2^*_{\mu}-1}_{\lambda_{\varepsilon},\xi_{\varepsilon}}, &  \text{  in  }\Omega_\varepsilon,\\
\omega=0, &  \text{  on  }\partial\Omega_\varepsilon.
\end{cases}
\end{align}
For simplicity, we will leave out the constant $\mathcal{A}_{H,L}$ in the sequel.
We look for a solution to problem (\ref{crihartree}) as
\begin{align}\label{vphi}
u_\varepsilon(x)=PU_{\lambda_{\varepsilon},\xi_{\varepsilon}}(x)
+\phi(x),
\end{align}
where the weight $\lambda_{\varepsilon}$ of the bubble and the center $\xi_{\varepsilon}$ of the bubble satisfy
\begin{align}\label{lambdaxijvarepsilon}
\lambda_{\varepsilon}=\varepsilon^{-\frac{1}{2}}\lambda \  \text{  and   } \ \xi_{\varepsilon}=\lambda^{-1}_{\varepsilon}\tau
\end{align}
for some positive number $\lambda$ and some point $\tau$ in $\mathbb{R}^N$. Moreover, we assume that given $\eta>0$ small
\begin{align}\label{lambdajtauj}
\eta<\lambda<\eta^{-1},\quad |\tau|<\eta.
\end{align}

Then $\phi$ satisfies
\begin{align}\label{phieq}
\begin{cases}
L_\varepsilon(\phi)=N_\varepsilon(\phi)+l_\varepsilon \ &\text{ in } \Omega_\varepsilon,\\
\phi=0\ &\text{ on } \partial\Omega_\varepsilon,
\end{cases}
\end{align}
where the linear operator $L_\varepsilon$ is defined by
\begin{align*}
L_\varepsilon(\phi)=-\Delta\phi-f'(PU_{\lambda_{\varepsilon},\xi_{\varepsilon}})\phi,
\end{align*}
the nonlinear part $N_\varepsilon$ is defined by
\begin{align*}
N_\varepsilon(\phi)=f(PU_{\lambda_{\varepsilon},\xi_{\varepsilon}}+\phi)-f(PU_{\lambda_{\varepsilon},\xi_{\varepsilon}})
-f'(PU_{\lambda_{\varepsilon},\xi_{\varepsilon}})\phi,
\end{align*}
the remainder term $l_\varepsilon$ is defined by
\begin{align*}
l_\varepsilon=-f(U_{\lambda_{\varepsilon},\xi_{\varepsilon}})+f(PU_{\lambda_{\varepsilon},\xi_{\varepsilon}}).
\end{align*}
Here $f(u)=\left(\int_{\Omega_\varepsilon}\frac{u^{2_{\mu}^*}}{|x-y|^{\mu}}dy\right)u^{2_{\mu}^*-1}$.

We denote by $PZ^h_{\lambda_{\varepsilon},\xi_{\varepsilon}}$ the projection of $Z^h_{\lambda_{\varepsilon},\xi_{\varepsilon}}$ onto $H^{1}_0(\Omega_\varepsilon)$ and define the subspace of $H^{1}_0(\Omega_\varepsilon)$
\begin{align}\label{K}
K:=\text{span}\left\{PZ^h_{\lambda_{\varepsilon},\xi_{\varepsilon}}:h=0,\cdots,N\right\},
\end{align}
and
\begin{align}\label{Kperp}
K^\perp:=\Bigg\{\phi\in H^{1}_0(\Omega_\varepsilon):\langle \phi,PZ^h_{\lambda_{\varepsilon},\xi_{\varepsilon}}\rangle=0,
h=0,\cdots,N\Bigg\}.
\end{align}

To prove our result we follow well known Liapunov-Schmidt reduction. First of all, we solve the following problem:
given a parameter $\lambda>0$ and a point $\tau=(\tau^1,\cdots,\tau^N)\in \mathbb{R}^N$, find a function $\phi:=\phi(\lambda, \tau)\in K^{\perp}$, such that
\begin{align}\label{linearprorl}
L_{\varepsilon}(\phi)=N_{\varepsilon}(\phi)+l_\varepsilon+\sum_{j=0}^Nc^jf'(PU_{\lambda_{\varepsilon},\xi_{\varepsilon}})
Z_{\lambda_{\varepsilon},\xi_{\varepsilon}}^j,
\end{align}
for certain constants $c^j$, depending on $\lambda$ and $\tau$, $j=0,\cdots, N$.

In fact we get the following existence result to the problem (\ref{linearprorl}).
\begin{proposition}\label{phi}
For any $\eta>0$ small but fixed, there exist $\varepsilon_0>0$ and $c>0$ such that for any $\lambda\in\mathbb{R}_{+}$ and $\tau\in\mathbb{R}^{N}$ satisfying (\ref{lambdajtauj}) and for any $\varepsilon \in(0,\varepsilon_0)$, there exists a unique solution $\phi=\phi(\lambda,\tau)$ which solves problem (\ref{linearprorl}). Moreover,
\begin{align*}
\|\phi\|=
\begin{cases}
O\left(\varepsilon^{\frac{N-\mu+2}{4}}\right), \ &\text{ if } N>6-\mu,\\
O\left(\varepsilon^{\frac{N-2}{2}}\left|\ln \varepsilon\right|\right), \ &\text{ if } N=6-\mu,\\
O\left(\varepsilon^{\frac{N-2}{2}}\right), \ &\text{ if } N<6-\mu.
\end{cases}
\end{align*}
\end{proposition}

Now, we can reduce the problem to a finite dimensional one. Define the functional
\begin{align*}
I_\varepsilon(u)
:=\frac{1}{2}\int_{\Omega_\varepsilon}|\nabla u|^2
-\frac{1}{2\cdot2_\mu^*}\int_{\Omega_\varepsilon}\int_{\Omega_\varepsilon}\frac{u^{2_\mu^*}(y)u^{2_\mu^*}(x)}{|x-y|^{\mu}}dxdy
\end{align*}
and
\begin{align*}
\tilde{I}_\varepsilon(\lambda,\tau):
=I_\varepsilon(PU_{\lambda_{\varepsilon},\xi_{\varepsilon}}+\phi).
\end{align*}

We denote by $G(x,y)$ the Green function of the Laplace operator in $\Omega$ with zero Dirichlet boundary condition, namely
\begin{align*}
\begin{cases}
-\Delta G(x,\cdot)=\delta(x-\cdot) &\text{ in }\Omega,\\
G(x,\cdot)=0 &\text{ on }\partial\Omega,
\end{cases}
\end{align*}
and we denote by $H(x,y)$ its regular part, i.e.,
\begin{align*}
H(x,\xi)=\frac{1}{(N-2)\omega_N|x-\xi|^{N-2}}-G(x,\xi),
\end{align*}
where $\omega_N$ is a measure of the unit sphere of $\mathbb{R}^N$ and $\delta(x)$ denotes the Dirac measure at the origin.
The Robin's function is defined as $H(x,x), x\in \Omega$.

One can get the following results.

\begin{proposition}\label{criticalpoint}
(i) If $(\bar{\lambda},\bar{\tau})$ is a critical point of $\tilde{I}_\varepsilon$, then $PU_{\lambda_{\varepsilon},\xi_{\varepsilon}}+\phi$ is a solution to problem (\ref{crihartree}).

(ii)
There holds
\begin{align*}
\tilde{I}_\varepsilon(\lambda,\tau)
=\left(1-\frac{1}{2_\mu^*}\right)\frac{N(N-2)}{2\mathcal{A}_{H,L}}A_N
+\frac{N(N-2)}{2\mathcal{A}_{H,L}}
\Psi({\tau},{\lambda})\varepsilon^{\frac{N-2}{2}}\left(1+o(1)\right)
\end{align*}
as $\varepsilon\to0$, $C^1$-uniformly with respect to $\lambda$ and $\tau$ satisfying (\ref{lambdajtauj}),
where
\begin{align*}
\Psi({\tau},{\lambda})
=\frac{(N-2)\omega_{N}B_NH(0,0)}{\lambda^{N-2}}
+\frac{M(\tau)}{(1+|\tau|^2)^{\frac{N-2}{2}}}
\lambda^{N-2}.
\end{align*}
and
\begin{equation}\label{e4.17}
A_N=\int_{\mathbb{R}^N}U_{1, 0}^{2^*}(x)dx, \quad B_N=\int_{\mathbb{R}^N}U_{1, 0}^{2^*-1}(x)dx, \quad
M(\tau)=\int_{\mathbb{R}^N}\frac{1}{|z|^{N-2}(1+|z-\tau|^2)^{\frac{N+2}{2}}}dz.
\end{equation}
\end{proposition}

Before proving Theorem \ref{thm1.1}, we should prove that $\tilde{I}_\varepsilon$ has a non degenerate critical point. To this end, we have the following proposition.
\begin{proposition}\label{nondectiticalpoint}
There exists $\bar{\lambda}_0\in \mathbb{R}^k_+$ such that $(0,\bar{\lambda}_0)$ is a non degenerate critical point of the function $\Psi(\tau,\lambda)$, which is defined by
\begin{align*}
\Psi(\tau,\lambda)
=\frac{(N-2)\omega_{N}B_NH(0,0)}{\lambda^{N-2}}
+\frac{M(\tau)}{(1+|\tau|^2)^{\frac{N-2}{2}}}
\lambda^{N-2}.
\end{align*}
\end{proposition}
\begin{proof}
Let us rewrite the function $\Psi(\tau,\lambda)$ as
\begin{align*}
\Psi^*(\tau,\mu)
:=m\mu^2+\frac{g(\tau)}{\mu^2},
\end{align*}
where $m=(N-2)\omega_{N}B_NH(0,0)$, $\mu=\frac{1}{\lambda^{\frac{N-2}{2}}}$, $g(\tau)=\frac{M(\tau)}{(1+|\tau|^2)^{\frac{N-2}{2}}}$.

In \cite[Lemma 4.1]{gemp} it is proved that
$\tau=0$ is a critical point for the function
$M(\tau)$. So, it
is also a critical point for $g(\tau)$. In addition,
If we fix $\tau=0$, the quadratic form $\mu\rightarrow m\mu^2$ is strictly positively definite, then the function $\mu\rightarrow \Psi^*(0,\mu)$ has a minimum point $\bar{\mu}_0$. Thus, we get that $(0,\bar{\mu}_0)$ is critical point for the function $\Psi^*$. Also, we can prove that it is non degenerate. In fact, we can write
\begin{align*}
\mathscr{H}\Psi^*(0,\bar{\mu}_0)=
\left(
  \begin{array}{cc}
    \mathscr{H}_{{\tau}}\Psi^*(0,\bar{\mu}_0) & 0 \\
    0 & \mathscr{H}_{{\mu}}\Psi^*(0,\bar{\mu}_0) \\
  \end{array}
\right),
\end{align*}
we get $\left|\mathscr{H}_{{\tau}}\Psi^*(0,\bar{\mu}_0)\right|\neq 0$. So, we are left to prove that $\left|\mathscr{H}_{{\mu}}\Psi^*(0,\bar{\mu}_0)\right|\neq 0$. In fact,
\begin{align*}
\partial_{\mu}\Psi^*(0,\bar{\mu}_0)
=2m\bar{\mu}_0-\frac{2g(0)}{\bar{\mu}_0^3}=0.
\end{align*}
Hence,
\begin{align*}
m\bar{\mu}_0^2=\frac{g(0)}{\bar{\mu}_0^2},
\end{align*}
and
\begin{align*}
\partial_{\mu\mu}\Psi^*(0,\bar{\mu}_0)
=2m+\frac{6g(0)}{\bar{\mu}_0^4}
=8m>0.
\end{align*}
Then, we deduce that
$\left|\mathscr{H}_{{\mu}}\Psi^*(0,\bar{\mu}_0)\right|\neq 0$. This implies that $\left|\mathscr{H}\Psi^*(0,\bar{\mu}_0)\right|\neq 0$. Then we conclude our proof.
\end{proof}

Now we are ready to prove our main result.

\noindent
{\it Proof of Theorem \ref{thm1.1}:} Combining Proposition \ref{phi} with $(i)$ in Proposition \ref{criticalpoint}, we see that if we want to prove that  $PU_{\lambda_{\varepsilon},\xi_{\varepsilon}}+\phi$ is the solution of problem (\ref{crihartree}), we need to find a critical point of the function $\tilde{I}_\varepsilon$.

By Lemma \ref{nondectiticalpoint}, there exists $\bar{\lambda}_0\in \mathbb{R}_+$
such that $(\bar{\lambda}_0,0)$ is a non degenerate critical point of the function $\Psi({\tau},{\lambda})$ defined in Proposition \ref{criticalpoint}, which is stable with respect to $C^1$-perturbation. Therefore, from $(ii)$ in Proposition \ref{criticalpoint}, we deduce that the function $\tilde{I}_\varepsilon$ has a critical point, denoted by $(\bar{\lambda},\bar{\tau})$, which satisfies $\bar{\lambda}\rightarrow \bar{\lambda}_0$ and $\bar{\tau}\rightarrow 0$ as $\varepsilon\rightarrow 0$.
\qed

\section{The Linear Problem}\label{The Linear Problem}

In this section, we mainly deal with the proof of Proposition \ref{phi}. To carry out the reduction argument, we define the projection $Q_\varepsilon$ from $H_0^1(\Omega_\varepsilon)$ to $K^{\perp}$ as follows
\begin{align*}
Q_{\varepsilon}u=u-\sum_{j=0}^Nc_{\varepsilon,j}PZ^j_{\lambda_{\varepsilon},\xi_{\varepsilon}}.
\end{align*}
Then we first study the invertibility of the operator $Q_{\varepsilon}L_\varepsilon$
and have the following lemma.
\begin{lemma}\label{Linver}
For any $\eta>0$ small but fixed, there exists $\varepsilon_0>0$ and $C>0$ such that for any ${\tau}\in \mathbb{R}^{N}$, ${\lambda}\in \mathbb{R}_+$ satisfying (\ref{lambdajtauj}) and for any $\varepsilon\in (0,\varepsilon_0)$, we have
\begin{align*}
\|Q_{\varepsilon} L_\varepsilon(\phi)\|\geq C\|\phi\| \quad \text{ for any } \phi \in K^\perp.
\end{align*}
\end{lemma}
\begin{proof}
By contradiction, suppose that there exist sequences $\varepsilon_n\rightarrow 0$, ${\tau}_n\in \mathbb{R}^{N}$, ${\lambda}_n\in \mathbb{R}_+$ where $\tau_{n}\rightarrow \tau\in \mathbb{R}^{N}$, with $|\tau|\leq \eta$, and $\lambda_{n}\rightarrow \lambda>0$, and $\phi_n\in K^{\perp}$, $\|\phi_n\|=1$ such that
\begin{align*}
\|Q_{\varepsilon_n} L_{\varepsilon_n}(\phi_n)\|\leq \frac{1}{n}\|\phi_n\|=\frac{1}{n}.
\end{align*}
For simplicity, we let $PU_{\lambda_{\varepsilon_n},\xi_{\varepsilon_n}}:=PU_n$.
Then for any $\varphi\in K^\perp$, we have
\begin{align}\label{linver2}
&\int_{\Omega_\varepsilon}\nabla \phi_n\nabla \varphi-2_{\mu}^*\int_{\Omega_\varepsilon}\int_{\Omega_\varepsilon}\frac{PU_n(y)^{2_{\mu}^*-1}\phi_n(y)
PU_n(x)^{2_{\mu}^*-1}\varphi(x)}
{|x-y|^{\mu}}dxdy
\nonumber\\
&-(2_{\mu}^*-1)\int_{\Omega_\varepsilon}\int_{\Omega_\varepsilon}\frac{PU_n(y)^{2_{\mu}^*}
PU_n(x)^{2_{\mu}^*-2}\phi_n(x)\varphi(x)}
{|x-y|^{\mu}}dxdy
\nonumber\\
=&\left\langle Q_{\varepsilon_n} L_{\varepsilon_n}(\phi_n), \varphi\right\rangle
\leq o(1)\|\varphi\|\|\phi_n\|=o(1)\|\varphi\|.
\end{align}
Choose $\varphi=\phi_n$ in (\ref{linver2}), then
\begin{align}\label{linver3}
&\int_{\Omega_\varepsilon}|\nabla \phi_n|^2-2_{\mu}^*\int_{\Omega_\varepsilon}\int_{\Omega_\varepsilon}\frac{PU_n(y)^{2_{\mu}^*-1}\phi_n(y)
PU_n(x)^{2_{\mu}^*-1}\phi_n(x)}
{|x-y|^{\mu}}dxdy
\nonumber\\
&-(2_{\mu}^*-1)\int_{\Omega_\varepsilon}\int_{\Omega_\varepsilon}\frac{PU_n(y)^{2_{\mu}^*}
PU_n(x)^{2_{\mu}^*-2}\phi_n(x)\phi_n(x)}
{|x-y|^{\mu}}dxdy
=o(1).
\end{align}
Next, we define
\begin{align*}
\tilde{\phi}_n(x)=\lambda_{\varepsilon_n}^{-\frac{N-2}{2}}
\phi_n(\lambda_{\varepsilon_n}^{-1}x+\xi_{\varepsilon_n}).
\end{align*}
We get $\int_{\Omega_\varepsilon}|\nabla \tilde{\phi}_n|^2\leq C$. Thus, up to a subsequence, we may assume that  $\tilde{\phi}_n\rightharpoonup \tilde{\phi}$ in $D^{1,2}(\mathbb{R}^N)$ and strongly in $L^q_{loc}(\mathbb{R}^N)$ for any $q\in [2,2^*)$. Our purpose is to prove that $\tilde{\phi}=0$, which is equivalent to prove the following problem
\begin{align}\label{linverlimiteq}
\begin{cases}
-\Delta \tilde{\phi}-2_{\mu}^*\Bigg(\frac{1}{|x|^{\mu}}\ast (U_{1,0}(y)^{2_{\mu}^*-1}\tilde{\phi})\Bigg)U_{1,0}(y)^{2_{\mu}^*-1}
-(2_{\mu}^*-1)\Bigg(\frac{1}{|x|^{\mu}}\ast U_{1,0}(y)^{2_{\mu}^*}\Bigg)U_{1,0}(y)^{2_{\mu}^*-2}\tilde{\phi}
=0,\\
\langle\frac{\partial U_{\lambda,0}}{\partial \lambda}\Bigg|_{\lambda=1},\tilde{\phi}\rangle
=\langle \frac{\partial U_{1,\xi}}{\partial \xi^j}\Bigg|_{\xi=0},\tilde{\phi}\rangle=0, \quad \text{ for } j=1,\cdots,N.
\end{cases}
\end{align}
We first prove the first identity in (\ref{linverlimiteq}).
For any $\varphi\in D^{1,2}(\mathbb{R}^N)$, we define
\begin{align*}
Q_{\varepsilon_n}\varphi=\varphi-\sum_{j=0}^Nc_{\varepsilon_n,j}PZ^j_{\lambda_{\varepsilon_n},\xi_{\varepsilon_n}},
\end{align*}
then we have
\begin{align*}
\begin{cases}
\sum_{j=0}^Nc_{\varepsilon_n,j}\left\langle PZ^j_{\lambda_{\varepsilon_n},\xi_{\varepsilon_n}}, PZ^i_{\lambda_{\varepsilon_n},\xi_{\varepsilon_n}}\right\rangle
=\left\langle \varphi, PZ^0_{\lambda_{\varepsilon_n},\xi_{\varepsilon_n}}\right\rangle,\\
\sum_{j=0}^Nc_{\varepsilon_n,j}\left\langle PZ^j_{\lambda_{\varepsilon_n},\xi_{\varepsilon_n}}, PZ^i_{\lambda_{\varepsilon_n},\xi_{\varepsilon_n}}\right\rangle
=\left\langle \varphi, PZ^0_{\lambda_{\varepsilon_n},\xi_{\varepsilon_n}}\right\rangle.
\end{cases}
\end{align*}
From the above identities, we see that
\begin{align*}
c_{\varepsilon_n,0}=a_{N}\left\langle \varphi, PZ^0_{\lambda_{\varepsilon_n},\xi_{\varepsilon_n}}\right\rangle,\quad
c_{\varepsilon_n,0}=b_{N,i}\left\langle \varphi, PZ^i_{\lambda_{\varepsilon_n},\xi_{\varepsilon_n}}\right\rangle
\end{align*}
for some constants $a_{N}$ and $b_{N,i}$. Therefore,
\begin{align}\label{linver4}
&\int_{\Omega_\varepsilon}\nabla \phi_n \nabla \varphi-2_{\mu}^*\int_{\Omega_\varepsilon}\int_{\Omega_\varepsilon}\frac{PU_n(y)^{2_{\mu}^*-1}\phi_n(y)
PU_n(x)^{2_{\mu}^*-1}\varphi(x)}
{|x-y|^{\mu}}dxdy
\nonumber\\
&-(2_{\mu}^*-1)\int_{\Omega_\varepsilon}\int_{\Omega_\varepsilon}\frac{PU_n(y)^{2_{\mu}^*}
PU_n(x)^{2_{\mu}^*-2}\phi_n(x)\varphi(x)}
{|x-y|^{\mu}}dxdy
\nonumber\\
=&\left\langle L_{\varepsilon_n}(\phi_n), \varphi\right\rangle
=\left\langle L_{\varepsilon_n}(\phi_n), Q_{\varepsilon_n}\varphi\right\rangle
+\sum_{j=0}^Nc_{\varepsilon_n,j}\left\langle L_{\varepsilon_n}(\phi_n), PZ^j_{\lambda_{\varepsilon_n},\xi_{\varepsilon_n}}\right\rangle
\nonumber\\
=&o(1)\|\varphi\|+\sum_{j=0}^Nc_{\varepsilon_n,j}\left\langle L_{\varepsilon_n}(\phi_n), PZ^j_{\lambda_{\varepsilon_n},\xi_{\varepsilon_n}}\right\rangle
\nonumber\\
=&o(1)\|\varphi\|+\tilde{a}_{N}\left\langle \varphi, PZ^0_{\lambda_{\varepsilon_n},\xi_{\varepsilon_n}}\right\rangle
+\sum_{i=1}^N\tilde{b}_{N,i}\left\langle \varphi, PZ^i_{\lambda_{\varepsilon_n},\xi_{\varepsilon_n}}\right\rangle.
\end{align}
Then taking $\varphi=PZ^i_{\lambda_{\varepsilon_n},\xi_{\varepsilon_n}}\in K^\perp$ in (\ref{linver4}) and using Lemma \ref{PU}, H\"{o}lder inequality and Hardy-Littlewood-Sobolev inequality, we have that for $i,j\neq 0$
\begin{align}\label{linver5}
&\sum_{j=1}^N\tilde{b}_{N,i}\left\langle PZ^j_{\lambda_{\varepsilon_n},\xi_{\varepsilon_n}}, PZ^i_{\lambda_{\varepsilon_n},\xi_{\varepsilon_n}}\right\rangle
\nonumber\\
=&\int_{\Omega_\varepsilon}\nabla \phi_n \nabla PZ^i_{\lambda_{\varepsilon_n},\xi_{\varepsilon_n}}
-2_{\mu}^*\int_{\Omega_\varepsilon}\int_{\Omega_\varepsilon}\frac{PU_n(y)^{2_{\mu}^*-1}\phi_n(y)
PU_n(x)^{2_{\mu}^*-1}PZ^i_{\lambda_{\varepsilon_n},\xi_{\varepsilon_n}}(x)}
{|x-y|^{\mu}}dxdy
\nonumber\\
&-(2_{\mu}^*-1)\int_{\Omega_\varepsilon}\int_{\Omega_\varepsilon}\frac{PU_n(y)^{2_{\mu}^*}
PU_n(x)^{2_{\mu}^*-2}\phi_n(x)PZ^i_{\lambda_{\varepsilon_n},\xi_{\varepsilon_n}}(x)}
{|x-y|^{\mu}}dxdy
-\tilde{a}_{N}\left\langle PZ^i_{\lambda_{\varepsilon_n},\xi_{\varepsilon_n}}, PZ^0_{\lambda_{\varepsilon_n},\xi_{\varepsilon_n}}\right\rangle
\nonumber\\
&+o(1)\|PZ^i_{\lambda_{\varepsilon_n},\xi_{\varepsilon_n}}\|
\nonumber\\
=&2_{\mu}^*\int_{\Omega_\varepsilon}\int_{\Omega_\varepsilon}\frac{PU_n(y)^{2_{\mu}^*-1}\phi_n(y)
PU_n(x)^{2_{\mu}^*-1}(Z^i_{\lambda_{\varepsilon_n},\xi_{\varepsilon_n}}-PZ^i_{\lambda_{\varepsilon_n},\xi_{\varepsilon_n}})(x)}
{|x-y|^{\mu}}dxdy
\nonumber\\
&+2_{\mu}^*\int_{\Omega_\varepsilon}\int_{\Omega_\varepsilon}\frac{\left(
U_n(y)^{2_{\mu}^*-1}\phi_n(y)U_n(x)^{2_{\mu}^*-1}
-PU_n(y)^{2_{\mu}^*-1}\phi_n(y)PU_n(x)^{2_{\mu}^*-1}
\right)
Z^i_{\lambda_{\varepsilon_n},\xi_{\varepsilon_n}}(x)}
{|x-y|^{\mu}}dxdy
\nonumber\\
&+(2_{\mu}^*-1)\int_{\Omega_\varepsilon}\int_{\Omega_\varepsilon}\frac{PU_n(y)^{2_{\mu}^*}
PU_n(x)^{2_{\mu}^*-2}\phi_n(x)(Z^i_{\lambda_{\varepsilon_n},\xi_{\varepsilon_n}}
-PZ^i_{\lambda_{\varepsilon_n},\xi_{\varepsilon_n}})(x)}
{|x-y|^{\mu}}dxdy
\nonumber\\
&+(2_{\mu}^*-1)\int_{\Omega_\varepsilon}\int_{\Omega_\varepsilon}\frac{
\left(U_n(y)^{2_{\mu}^*}U_n(x)^{2_{\mu}^*-2}\phi_n(x)
-PU_n(y)^{2_{\mu}^*}PU_n(x)^{2_{\mu}^*-2}\phi_n(x)
\right)
Z^i_{\lambda_{\varepsilon_n},\xi_{\varepsilon_n}}(x)}
{|x-y|^{\mu}}dxdy
\nonumber\\
&+o(1)\|PZ^i_{\lambda_{\varepsilon_n},\xi_{\varepsilon_n}}\|
\nonumber\\
=&
O\left(\left(\int_{\Omega_\varepsilon}\left|Z^i_{\lambda_{\varepsilon_n},\xi_{\varepsilon_n}}
-PZ^i_{\lambda_{\varepsilon_n},\xi_{\varepsilon_n}}
\right|^{\frac{2N}{N-2}}\right)^{\frac{N-2}{2N}}\right)
+O\left(\left(\int_{\Omega_\varepsilon}\left|U_n^{2_{\mu}^*-1}-PU_n^{2_{\mu}^*-1}
\right|^{\frac{2N}{N-\mu+2}}\right)^{\frac{N-\mu+2}{2N}}\right)\|PZ^i_{\lambda_{\varepsilon_n},\xi_{\varepsilon_n}}\|
\nonumber\\
&+O\left(\left(\int_{\Omega_\varepsilon}\left|U_n^{2_{\mu}^*}-PU_n^{2_{\mu}^*}
\right|^{\frac{2N}{2N-\mu}}\right)^{\frac{2N-\mu}{2N}}\right)\|PZ^i_{\lambda_{\varepsilon_n},\xi_{\varepsilon_n}}\|
+O\left(\left(\int_{\Omega_\varepsilon}\left|U_n^{2_{\mu}^*-2}-PU_n^{2_{\mu}^*-2}
\right|^{\frac{2N}{4-\mu}}\right)^{\frac{4-\mu}{2N}}\right)\|PZ^i_{\lambda_{\varepsilon_n},\xi_{\varepsilon_n}}\|
\nonumber\\
&+o(1)\|PZ^i_{\lambda_{\varepsilon_n},\xi_{\varepsilon_n}}\|
\nonumber\\
=&O\left(\frac{1}{\lambda_{\varepsilon_n}^{\frac{N-2}{2}}}\right)\|PZ^i_{\lambda_{\varepsilon_n},\xi_{\varepsilon_n}}\|+o\left(\lambda_{\varepsilon_n}\right).
\end{align}
Hence,
\begin{align*}
\tilde{b}_{N,i}=O\left(\frac{1}{\lambda_{\varepsilon_n}^{\frac{N-4}{2}}}\right)+o\left(\frac{1}{\lambda_{\varepsilon_n}}\right),
\end{align*}
and
\begin{align}\label{an}
\tilde{b}_{N,i}\left\langle \varphi, PZ^i_{\lambda_{\varepsilon_n},\xi_{\varepsilon_n}}\right\rangle
=\tilde{b}_{N,i}O\left(\|PZ^i_{\lambda_{\varepsilon_n},\xi_{\varepsilon_n}}\|\|\varphi\|\right)
=o(1)\|\varphi\|.
\end{align}
Similarly, if we choose $\varphi=PZ^0_{\lambda_{\varepsilon_n},\xi_{\varepsilon_n}}$, we can get
\begin{align*}
\tilde{a}_{N}=O\left(\frac{1}{\lambda_{\varepsilon_n}^{\frac{N}{2}}}\right)+o\left(\lambda_{\varepsilon_n}\right),
\end{align*}
\begin{align}\label{bn}
\tilde{a}_{N}\left\langle \varphi, PZ^0_{\lambda_{\varepsilon_n},\xi_{\varepsilon_n}}\right\rangle
=o(1)\|\varphi\|.
\end{align}
In terms of (\ref{linver4}), (\ref{an}), (\ref{bn}), we get
\begin{align}\label{linver4limit}
&\int_{\Omega_\varepsilon}\nabla \phi_n \nabla \varphi-2_{\mu}^*\int_{\Omega_\varepsilon}\int_{\Omega_\varepsilon}\frac{PU_n(y)^{2_{\mu}^*-1}\phi_n(y)
PU_n(x)^{2_{\mu}^*-1}\varphi(x)}
{|x-y|^{\mu}}dxdy
\nonumber\\
&-(2_{\mu}^*-1)\int_{\Omega_\varepsilon}\int_{\Omega_\varepsilon}\frac{PU_n(y)^{2_{\mu}^*}
PU_n(x)^{2_{\mu}^*-2}\phi_n(x)\varphi(x)}
{|x-y|^{\mu}}dxdy
\nonumber\\
=&o(1)\|\varphi\|.
\end{align}
If we set $\tilde{\varphi}(x)=\lambda_{\varepsilon_n}^{\frac{N-2}{2}}\varphi(\lambda_{\varepsilon_n}(x-\xi_{\varepsilon_n}))$, then by (\ref{linver4limit}), we obtain
\begin{align*}
&\int_{\lambda_{\varepsilon_n}(\Omega_\varepsilon-\xi_{\varepsilon_n})}\nabla \tilde{\phi}_n \nabla \varphi-2_{\mu}^*\int_{\lambda_{\varepsilon_n}(\Omega_\varepsilon-\xi_{\varepsilon_n})}
\int_{\lambda_{\varepsilon_n}(\Omega_\varepsilon-\xi_{\varepsilon_n})}\frac{U_{0,1}(y)^{2_{\mu}^*-1}\tilde{\phi}_n(y)
U_{0,1}(x)^{2_{\mu}^*-1}\varphi(x)}
{|x-y|^{\mu}}dxdy
\nonumber\\
&-(2_{\mu}^*-1)\int_{\lambda_{\varepsilon_n}(\Omega_\varepsilon-\xi_{\varepsilon_n})}
\int_{\lambda_{\varepsilon_n}(\Omega_\varepsilon-\xi_{\varepsilon_n})}\frac{U_{0,1}(y)^{2_{\mu}^*}
U_{0,1}(x)^{2_{\mu}^*-2}\tilde{\phi}_n(x)\varphi(x)}
{|x-y|^{\mu}}dxdy
\nonumber\\
=&\int_{\Omega_\varepsilon}\nabla \phi_n \nabla \tilde{\varphi}-2_{\mu}^*\int_{\Omega_\varepsilon}\int_{\Omega_\varepsilon}\frac{U_n(y)^{2_{\mu}^*-1}\phi_n(y)
U_n(x)^{2_{\mu}^*-1}\tilde{\varphi}(x)}
{|x-y|^{\mu}}dxdy
\nonumber\\
&-(2_{\mu}^*-1)\int_{\Omega_\varepsilon}\int_{\Omega_\varepsilon}\frac{U_n(y)^{2_{\mu}^*}
U_n(x)^{2_{\mu}^*-2}\phi_n(x)\tilde{\varphi}(x)}
{|x-y|^{\mu}}dxdy
\nonumber\\
=&o(1)\|\tilde{\varphi}\|=o(1).
\end{align*}
Taking $n\rightarrow \infty$, we know that
\begin{align*}
-\Delta \tilde{\phi}-2_{\mu}^*\Bigg(\frac{1}{|x|^{\mu}}\ast (U_{1,0}(y)^{2_{\mu}^*-1}\tilde{\phi})\Bigg)U_{1,0}(y)^{2_{\mu}^*-1}
-(2_{\mu}^*-1)\Bigg(\frac{1}{|x|^{\mu}}\ast U_{1,0}(y)^{2_{\mu}^*}\Bigg)U_{1,0}(y)^{2_{\mu}^*-2}\tilde{\phi}
=0.
\end{align*}

To prove the second identity in (\ref{linverlimiteq}) holds, it can follows from $\tilde{\phi}_n\rightharpoonup \tilde{\phi}$ and $\tilde{\phi}_n\in K^\perp$.
Therefore, we have $\tilde{\phi}=0$ by (\ref{linverlimiteq}). Moreover, by Hardy-Littlewood-Sobolev inequality and (\ref{linver3}), we have
\begin{align*}
\int_{\Omega_\varepsilon}|\nabla \tilde{\phi}_n|^2
=&2_{\mu}^*\int_{\Omega_\varepsilon}\int_{\Omega_\varepsilon}\frac{PU_n(y)^{2_{\mu}^*-1}\tilde{\phi}_n(y)
PU_n(x)^{2_{\mu}^*-1}\tilde{\phi}_n(x)}
{|x-y|^{\mu}}dxdy
\nonumber\\
&-(2_{\mu}^*-1)\int_{\Omega_\varepsilon}\int_{\Omega_\varepsilon}\frac{PU_n(y)^{2_{\mu}^*}
PU_n(x)^{2_{\mu}^*-2}\tilde{\phi}_n(x)\tilde{\phi}_n(x)}
{|x-y|^{\mu}}dxdy+o(1)
\nonumber\\
=&O\left(\left(\int_{\Omega_\varepsilon}\left(U_n(y)^{2_{\mu}^*-1}\tilde{\phi}_n\right)^{\frac{2N}{2N-\mu}}\right)^{\frac{2N-\mu}{N}}
+\left(\int_{\Omega_\varepsilon}\left(U_n(y)^{2_{\mu}^*-2}\tilde{\phi}_n^2\right)^{\frac{2N}{2N-\mu}}\right)^{\frac{2N-\mu}{2N}}
\right)
=o(1).
\end{align*}
Then we prove that
\begin{align*}
\int_{\Omega_\varepsilon}|\nabla \tilde{\phi}_n|^2
=o(1),
\end{align*}
which contradicts the assumption that
$\|\phi_n\|=1$. That conclude the proof.
\end{proof}

We are now ready to carry out the reduction for (\ref{phieq}). That is, we consider the following problem
\begin{align}\label{phieqQ}
\begin{cases}
Q_\varepsilon L_\varepsilon(\phi)=Q_\varepsilon N_\varepsilon(\phi)+Q_\varepsilon l_\varepsilon \ &\text{ in } \Omega_\varepsilon,\\
\phi=0\ &\text{ on } \partial\Omega_\varepsilon,
\end{cases}
\end{align}

In order to prove Proposition \ref{phi}, we have to estimate $l_\varepsilon$ and $R_\varepsilon$ respectively. To this end, the first order approximation of the function $PU_{\lambda,\xi}$ is given in the following, which is analogous to \cite[Lemma 3.1]{gemp}.

\begin{lemma}\label{PU}
Suppose that $\xi=\lambda^{-1}\tau$, with $\lambda\rightarrow \infty$ as $\varepsilon\rightarrow 0$ and $\varepsilon=o(\lambda^{-1})$ as $\varepsilon\rightarrow 0$. Then, if we define
\begin{align*}
R(x)=PU_{\lambda,\xi}-U_{\lambda,\xi}+\frac{(N-2)\omega_{N}}{\lambda^{\frac{N-2}{2}}}H(x,\xi)
+\frac{\lambda^{\frac{N-2}{2}}}{(1+|\tau|^2)^{\frac{N-2}{2}}}\frac{\varepsilon^{N-2}}{|x|^{N-2}},
\end{align*}
there exists a positive constant $c$ such that for any $x\in \Omega_\varepsilon$
\begin{align}\label{PU1}
|R(x)|\leq c\lambda^{-\frac{N-2}{2}}\Bigg[\frac{\varepsilon^{N-2}(1+\varepsilon\lambda^{N-1})}{|x|^{N-2}}
+\lambda^{-2}+(\varepsilon\lambda)^{N-2}
\Bigg],
\end{align}
\begin{align}\label{PU2}
|R_{\lambda}(x)|\leq c\lambda^{-\frac{N-4}{2}}\Bigg[\frac{\varepsilon^{N-2}(1+\varepsilon\lambda^{N-3})}{|x|^{N-2}}
+\lambda^{-2}\left(\lambda^{-2}+(\varepsilon\lambda)^{N-2}\right)
\Bigg],
\end{align}
\begin{align}\label{PU3}
|R_{\tau^i}(x)|\leq c\lambda^{-\frac{N}{2}}\Bigg[\frac{\varepsilon^{N-2}(1+\varepsilon\lambda^{N})}{|x|^{N-2}}
+\lambda^{-2}\left(\lambda^{-2}+\varepsilon^{N-1}\lambda^{N-2}\right)
\Bigg].
\end{align}
\end{lemma}
\begin{proof}
The proof is similar to \cite[Lemma 3.1]{gemp}, so we limit ourselves to prove (\ref{PU1}), the
other two claims could be obtained analogously.

Let us set $\hat{R}(x)=\lambda^{\frac{N-2}{2}}R(\varepsilon x)$, $x\in \left(\varepsilon^{-1}\Omega\backslash B(0,1)\right):=\hat{\Omega}_\varepsilon$ and $\hat{\Omega}_\varepsilon \rightarrow \mathbb{R}^N\backslash B(0,1)$ as $\varepsilon \rightarrow 0$. Then
\begin{align*}
\begin{cases}
-\Delta \hat{R}=0 & \text{  in } \hat{\Omega}_\varepsilon, \\
\hat{R}(x)=-\frac{\lambda^{N-2}}{(1+|\varepsilon \lambda x-\tau|^2)^{\frac{N-2}{2}}}
+(N-2)\omega_{N}H(\varepsilon x,\xi)
+\frac{\lambda^{N-2}}{(1+|\tau|^2)^{\frac{N-2}{2}}} & \text{  on } \partial B(0,1), \\
\hat{R}(x)=-\frac{\lambda^{N-2}}{(1+|\varepsilon \lambda x-\tau|^2)^{\frac{N-2}{2}}}
+\frac{1}{|\varepsilon x-\xi|^{N-2}}
+\frac{\lambda^{N-2}}{(1+|\tau|^2)^{\frac{N-2}{2}}|x|^{N-2}} & \text{  on } \partial(\varepsilon^{-1}\Omega).
\end{cases}
\end{align*}
Moreover, the following estimates hold true
\begin{align*}
|\hat{R}(x)|\leq c(1+\lambda^{N-2}\cdot\varepsilon \lambda)\quad \forall x\in \partial B(0,1),
\end{align*}
and
\begin{align*}
|\hat{R}(x)|\leq c(\lambda^{-2}+(\varepsilon \lambda)^{N-2})\quad \forall x\in \partial(\varepsilon^{-1}\Omega).
\end{align*}
By using a comparison argument, we infer
\begin{align}\label{hatR}
|\hat{R}(x)|\leq c\Bigg[\frac{(1+\varepsilon\lambda^{N-1})}{|x|^{N-2}}
+\lambda^{-2}+(\varepsilon\lambda)^{N-2}\Bigg] \quad \forall x\in \hat{\Omega}_\varepsilon.
\end{align}
Therefore, by (\ref{hatR}) we deduce (\ref{PU1}).

\end{proof}

\begin{lemma}\label{lvarepsilon}
For $N\geq 5$, it holds
\begin{align*}
\|l_\varepsilon\|=
\begin{cases}
O\left(\varepsilon^{\frac{N-\mu+2}{4}}\right), \ &\text{ if } N>6-\mu,\\
O\left(\varepsilon^{\frac{N-2}{2}}\left|\ln \varepsilon\right|\right), \ &\text{ if } N=6-\mu,\\
O\left(\varepsilon^{\frac{N-2}{2}}\right), \ &\text{ if } N<6-\mu.
\end{cases}
\end{align*}
\end{lemma}
\begin{proof}
One has
\begin{align*}
&\left\langle l_\varepsilon, \varphi\right\rangle
\nonumber\\
=&\int_{\Omega_\varepsilon}\left(-f(U_{\lambda_\varepsilon,\xi_\varepsilon})+f(PU_{\lambda_\varepsilon,\xi_\varepsilon})\right)\varphi
\nonumber\\
=&\int_{\Omega_\varepsilon}\int_{\Omega_\varepsilon}\frac{PU_{\lambda_\varepsilon,\xi_\varepsilon}^{2_{\mu}^*}(y)
PU_{\lambda_\varepsilon,\xi_\varepsilon}^{2_{\mu}^*-1}(x)\varphi(x)}{|x-y|^{\mu}}dxdy
-\int_{\Omega_\varepsilon}\int_{\Omega_\varepsilon}\frac{
U_{\lambda_\varepsilon,\xi_\varepsilon}^{2_{\mu}^*}(y)U_{\lambda_\varepsilon,\xi_\varepsilon}^{2_{\mu}^*-1}(x)\varphi(x)}{|x-y|^{\mu}}dxdy
\nonumber\\
=&\int_{\Omega_\varepsilon}\int_{\Omega_\varepsilon}\frac{\Bigg[PU_{\lambda_\varepsilon,\xi_\varepsilon}^{2_{\mu}^*}
-U_{\lambda_\varepsilon,\xi_\varepsilon}^{2_{\mu}^*}
\Bigg](y)PU_{\lambda_\varepsilon,\xi_\varepsilon}^{2_{\mu}^*-1}(x)\varphi(x)}{|x-y|^{\mu}}dxdy
\nonumber\\
&-\int_{\Omega_\varepsilon}\int_{\Omega_\varepsilon}\frac{U_{\lambda_\varepsilon,\xi_\varepsilon}^{2_{\mu}^*}(y)
\Bigg[PU_{\lambda_\varepsilon,\xi_\varepsilon}^{2_{\mu}^*-1}-U_{\lambda_\varepsilon,\xi_\varepsilon}^{2_{\mu}^*-1}
\Bigg](x)\varphi(x)}{|x-y|^{\mu}}dxdy
\nonumber\\
:=&l_1+l_2.
\end{align*}
By Hardy-Littlewood-Sobolev inequality, H\"{o}lder inequality, Sobolev inequality and \cite[Lemma 2.2]{yyli}, we get
\begin{align*}
\left|l_1\right|
\leq &C\left(\int_{\Omega_\varepsilon}\left|PU_{\lambda_\varepsilon,\xi_\varepsilon}^{2_{\mu}^*}-U_{\lambda_\varepsilon,\xi_\varepsilon}^{2_{\mu}^*}
\right|^{\frac{2N}{2N-\mu}}\right)^{\frac{2N-\mu}{2N}}
\left(\int_{\Omega_\varepsilon}\left|
PU_{\lambda_\varepsilon,\xi_\varepsilon}^{2_{\mu}^*-1}\varphi
\right|^{\frac{2N}{2N-\mu}}\right)^{\frac{2N-\mu}{2N}}
\nonumber\\
\leq &C
\left(\int_{\Omega_\varepsilon}\left|PU_{\lambda_\varepsilon,\xi_\varepsilon}^{2_{\mu}^*}-U_{\lambda_\varepsilon,\xi_\varepsilon}^{2_{\mu}^*}
\right|^{\frac{2N}{2N-\mu}}\right)^{\frac{2N-\mu}{2N}}
\left(\int_{\Omega_\varepsilon}\left|
PU_{\lambda_\varepsilon,\xi_\varepsilon}
\right|^{\frac{2N}{N-2}}\right)^{\frac{N-\mu+2}{2N}}
\left(\int_{\Omega_\varepsilon}\left|\varphi
\right|^{\frac{2N}{N-2}}\right)^{\frac{N-2}{2N}}
\nonumber\\
\leq &C
\left(\int_{\Omega_\varepsilon}\left|PU_{\lambda_\varepsilon,\xi_\varepsilon}^{2_{\mu}^*}-U_{\lambda_\varepsilon,\xi_\varepsilon}^{2_{\mu}^*}
\right|^{\frac{2N}{2N-\mu}}\right)^{\frac{2N-\mu}{2N}}\|\varphi\|
\nonumber\\
\leq &C\left(\int_{\Omega_\varepsilon}\left[U_{\lambda_\varepsilon,\xi_\varepsilon}^{(2_{\mu}^*-1)\frac{2N}{2N-\mu}}
\left|PU_{\lambda_\varepsilon,\xi_\varepsilon}-U_{\lambda_\varepsilon,\xi_\varepsilon}\right|^{\frac{2N}{2N-\mu}}
+\left|PU_{\lambda_\varepsilon,\xi_\varepsilon}-U_{\lambda_\varepsilon,\xi_\varepsilon}\right|^{2_{\mu}^*\frac{2N}{2N-\mu}}
\right]\right)^{\frac{2N-\mu}{2N}}\|\varphi\|,
\end{align*}
where it follows from Lemma \ref{PU} that
\begin{align*}
&\int_{\Omega_\varepsilon}\left[U_{\lambda_\varepsilon,\xi_\varepsilon}^{(2_{\mu}^*-1)\frac{2N}{2N-\mu}}
\left|PU_{\lambda_\varepsilon,\xi_\varepsilon}-U_{\lambda_\varepsilon,\xi_\varepsilon}\right|^{\frac{2N}{2N-\mu}}
\right]
\nonumber\\
=&\int_{\Omega_\varepsilon}\left[U_{\lambda_\varepsilon,\xi_\varepsilon}^{(2_{\mu}^*-1)\frac{2N}{2N-\mu}}
\left|\frac{(N-2)\omega_{N}}{\lambda_\varepsilon^{\frac{N-2}{2}}}H(x,\xi_\varepsilon)
+\frac{\lambda_\varepsilon^{\frac{N-2}{2}}}{(1+|\tau|^2)^{\frac{N-2}{2}}}\frac{\varepsilon^{N-2}}{|x|^{N-2}}
-R(x)
\right|^{\frac{2N}{2N-\mu}}
\right]
\nonumber\\
=&\int_{\Omega_\varepsilon}\left[
\left(\frac{\lambda_\varepsilon}
{1+\lambda_\varepsilon|x-\xi_\varepsilon|^2}
\right)^{\frac{N-2}{2}(2_{\mu}^*-1)\frac{2N}{2N-\mu}}
\left|\frac{(N-2)\omega_{N}}{\lambda_\varepsilon^{\frac{N-2}{2}}}H(x,\xi_\varepsilon)
+\frac{\lambda_\varepsilon^{\frac{N-2}{2}}}{(1+|\tau|^2)^{\frac{N-2}{2}}}\frac{\varepsilon^{N-2}}{|x|^{N-2}}
-R(x)
\right|^{\frac{2N}{2N-\mu}}
\right]
\nonumber\\
=&O\left(\frac{1}{\lambda_\varepsilon^{\frac{2N(N-2)}{2N-\mu}}}
\right)
+O\left(\left(\varepsilon\lambda_\varepsilon
\right)^{\frac{2N(N-2)}{2N-\mu}}\right)
=O\left(\varepsilon^{\frac{N(N-2)}{2N-\mu}}\right)
\end{align*}
and
\begin{align*}
&\int_{\Omega_\varepsilon}\left|PU_{\lambda_\varepsilon,\xi_\varepsilon}-U_{\lambda_\varepsilon,\xi_\varepsilon}\right|^{2_{\mu}^*\frac{2N}{2N-\mu}}
\nonumber\\
=&\int_{\Omega_\varepsilon}\left|\frac{(N-2)\omega_{N}}{\lambda_\varepsilon^{\frac{N-2}{2}}}H(x,\xi_\varepsilon)
+\frac{\lambda_\varepsilon^{\frac{N-2}{2}}}{(1+|\tau|^2)^{\frac{N-2}{2}}}\frac{\varepsilon^{N-2}}{|x|^{N-2}}
-R(x)
\right|^{2_{\mu}^*\frac{2N}{2N-\mu}}
\nonumber\\
=&O\left(\frac{1}{\lambda_\varepsilon^{N}}\right)
+O\left(\varepsilon^{2N}\lambda_\varepsilon^{N}\right)
=O\left(\varepsilon^{\frac{N}{2}}\right).
\end{align*}
Thus
\begin{align*}
\left|l_1\right|
=O\left(\frac{1}{\lambda_\varepsilon^{N-2}}
\right)
+O\left(\varepsilon^{N-2}\lambda_\varepsilon^{N-2}
\right)
=O\left(\varepsilon^{\frac{N-2}{2}}\right).
\end{align*}
In the same way, we also obtain
\begin{align*}
\left|l_2\right|\leq &C
\left(\int_{\Omega_\varepsilon}\left|PU_{\lambda_\varepsilon,\xi_\varepsilon}^{2_{\mu}^*-1}
-U_{\lambda_\varepsilon,\xi_\varepsilon}^{2_{\mu}^*-1}
\right|^{\frac{2N}{N-\mu+2}}\right)^{\frac{N-\mu+2}{2N}}\|\varphi\|
\nonumber\\
\leq &C
\left(\int_{\Omega_\varepsilon}U_{\lambda_\varepsilon,\xi_\varepsilon}^{(2_{\mu}^*-2)\frac{2N}{N-\mu+2}}
\left|PU_{\lambda_\varepsilon,\xi_\varepsilon}-U_{\lambda_\varepsilon,\xi_\varepsilon}\right|^{\frac{2N}{N-\mu+2}}
+\left|PU_{\lambda_\varepsilon,\xi_\varepsilon}-U_{\lambda_\varepsilon,\xi_\varepsilon}\right|^{(2_{\mu}^*-1)\frac{2N}{N-\mu+2}}
\right)^{\frac{N-\mu+2}{2N}}\|\varphi\|.
\end{align*}
Moreover, in terms of Lemma \ref{PU}, we have
\begin{align*}
&\int_{\Omega_\varepsilon}U_{\lambda_\varepsilon,\xi_\varepsilon}^{(2_{\mu}^*-2)\frac{2N}{N-\mu+2}}
\left|PU_{\lambda_\varepsilon,\xi_\varepsilon}-U_{\lambda_\varepsilon,\xi_\varepsilon}\right|
^{\frac{2N}{N-\mu+2}}
\nonumber\\
=&\int_{\Omega_\varepsilon}U_{\lambda_\varepsilon,\xi_\varepsilon}^{(2_{\mu}^*-2)\frac{2N}{N-\mu+2}}
\left|\frac{(N-2)\omega_{N}}{\lambda_\varepsilon^{\frac{N-2}{2}}}H(x,\xi_\varepsilon)
+\frac{\lambda_\varepsilon^{\frac{N-2}{2}}}{(1+|\tau|^2)^{\frac{N-2}{2}}}\frac{\varepsilon^{N-2}}{|x|^{N-2}}
-R(x)
\right|^{\frac{2N}{N-\mu+2}}
\nonumber\\
=&O\Bigg(\Bigg(\varepsilon \lambda_\varepsilon \Bigg)^{\frac{2N(N-2)}{N-\mu+2}}\Bigg)+
\begin{cases}
O\Bigg(\frac{1}{\lambda_\varepsilon^{\frac{2N(N-2)}{N-\mu+2}}}
\Bigg), \ &\text{ if } N<6-\mu,\\
O\Bigg(\frac{1}{\lambda_\varepsilon^{\frac{2N(N-2)}{N-\mu+2}}}
|\ln \lambda_\varepsilon|^{\frac{2N}{N-\mu+2}}\Bigg), \ &\text{ if } N=6-\mu,\\
O\Bigg(\frac{1}{\lambda_\varepsilon^{N}}\Bigg), \ &\text{ if } N>6-\mu,
\end{cases}
\nonumber\\
=&O\Bigg(\varepsilon^{\frac{N(N-2)}{N-\mu+2}}\Bigg)+
\begin{cases}
O\Bigg(\varepsilon^{\frac{N(N-2)}{N-\mu+2}}
\Bigg), \ &\text{ if } N<6-\mu,\\
O\Bigg(\varepsilon^{\frac{N(N-2)}{N-\mu+2}}
|\ln \varepsilon|^{\frac{2N}{N-\mu+2}}\Bigg), \ &\text{ if } N=6-\mu,\\
O\Bigg(\varepsilon^{\frac{N}{2}}\Bigg), \ &\text{ if } N>6-\mu,
\end{cases}
\end{align*}
and
\begin{align*}
&\int_{\Omega_\varepsilon}\left|PU_{\lambda_\varepsilon,\xi_\varepsilon}-U_{\lambda_\varepsilon,\xi_\varepsilon}\right|
^{(2_{\mu}^*-1)\frac{2N}{N-\mu+2}}
\nonumber\\
=&\int_{\Omega_\varepsilon}\left|
\frac{(N-2)\omega_{N}}{\lambda_\varepsilon^{\frac{N-2}{2}}}H(x,\xi_\varepsilon)
+\frac{\lambda_\varepsilon^{\frac{N-2}{2}}}{(1+|\tau|^2)^{\frac{N-2}{2}}}\frac{\varepsilon^{N-2}}{|x|^{N-2}}
-R(x)\right|^{(2_{\mu}^*-1)\frac{2N}{N-\mu+2}}
\nonumber\\
=&O\Bigg(\frac{1}{\lambda_\varepsilon^{N}}\Bigg)
+O\Bigg(\Bigg(\varepsilon \lambda_\varepsilon \Bigg)^{N}\varepsilon^{N}\Bigg)
=O\Bigg(\varepsilon^{\frac{N}{2}}\Bigg).
\end{align*}
Consequently,
\begin{align*}
\|l_\varepsilon\|=
\begin{cases}
O\left(\varepsilon^{\frac{N-\mu+2}{4}}\right), \ &\text{ if } N>6-\mu,\\
O\left(\varepsilon^{\frac{N-2}{2}}\left|\ln \varepsilon\right|\right), \ &\text{ if } N=6-\mu,\\
O\left(\varepsilon^{\frac{N-2}{2}}\right), \ &\text{ if } N<6-\mu.
\end{cases}
\end{align*}
\end{proof}

\begin{lemma}\label{R}
We have
\begin{align*}
\|N_\varepsilon(\phi)\|
\leq C\|\phi\|^{\min\{2,2_{\mu}^*-1\}}.
\end{align*}
\end{lemma}
\begin{proof}
Using \cite[Lemma A.1]{mussop} and the definition of $N_\varepsilon(\phi)$, we have that for $\varphi\in H_0^1(\mathbb{R}^N)$,
\begin{align*}
\left\langle N_\varepsilon(\phi),\varphi\right\rangle
=&\int_{\Omega_\varepsilon}\Bigg\{f(PU_{\lambda_\varepsilon,\xi_\varepsilon}+\phi)-f(PU_{\lambda_\varepsilon,\xi_\varepsilon})
-f'(PU_{\lambda_\varepsilon,\xi_\varepsilon})\phi\Bigg\}\varphi
\nonumber\\
=&\int_{\Omega_\varepsilon}\int_{\Omega_\varepsilon}\frac{(PU_{\lambda_\varepsilon,\xi_\varepsilon}+\phi)^{2_{\mu}^*}(y)
(PU_{\lambda_\varepsilon,\xi_\varepsilon}+\phi)^{2_{\mu}^*-1}(x)\varphi(x)}{|x-y|^{\mu}}dxdy
\nonumber\\
&-\int_{\Omega_\varepsilon}\int_{\Omega_\varepsilon}\frac{PU_{\lambda_\varepsilon,\xi_\varepsilon}^{2_{\mu}^*}(y)
PU_{\lambda_\varepsilon,\xi_\varepsilon}^{2_{\mu}^*-1}(x)\varphi(x)}{|x-y|^{\mu}}dxdy
\nonumber\\
&-2_{\mu}^*\int_{\Omega_\varepsilon}\int_{\Omega_\varepsilon}\frac{PU_{\lambda_\varepsilon,\xi_\varepsilon}^{2_{\mu}^*-1}(y)\phi(y)
PU_{\lambda_\varepsilon,\xi_\varepsilon}^{2_{\mu}^*-1}(x)\varphi(x)}{|x-y|^{\mu}}dxdy
\nonumber\\
&-(2_{\mu}^*-1)\int_{\Omega_\varepsilon}\int_{\Omega_\varepsilon}\frac{PU_{\lambda_\varepsilon,\xi_\varepsilon}^{2_{\mu}^*}(y)
PU_{\lambda_\varepsilon,\xi_\varepsilon}^{2_{\mu}^*-2}(x)\phi(x)\varphi(x)}{|x-y|^{\mu}}dxdy
\nonumber\\
=&\int_{\Omega_\varepsilon}\int_{\Omega_\varepsilon}\frac{\Bigg[(PU_{\lambda_\varepsilon,\xi_\varepsilon}+\phi)^{2_{\mu}^*}
-PU_{\lambda_\varepsilon,\xi_\varepsilon}^{2_{\mu}^*}
-2_{\mu}^*PU_{\lambda_\varepsilon,\xi_\varepsilon}^{2_{\mu}^*-1}\phi
\Bigg](y)
(PU_{\lambda_\varepsilon,\xi_\varepsilon}+\phi)^{2_{\mu}^*-1}(x)\varphi(x)}{|x-y|^{\mu}}dxdy
\nonumber\\
&+\int_{\Omega_\varepsilon}\int_{\Omega_\varepsilon}\frac{PU_{\lambda_\varepsilon,\xi_\varepsilon}^{2_{\mu}^*}(y)
\Bigg[(PU_{\lambda_\varepsilon,\xi_\varepsilon}+\phi)^{2_{\mu}^*-1}
-PU_{\lambda_\varepsilon,\xi_\varepsilon}^{2_{\mu}^*-1}
-(2_{\mu}^*-1)PU_{\lambda_\varepsilon,\xi_\varepsilon}^{2_{\mu}^*-2}\phi
\Bigg](x)\varphi(x)}{|x-y|^{\mu}}dxdy
\nonumber\\
&+2_{\mu}^*\int_{\Omega_\varepsilon}\int_{\Omega_\varepsilon}\frac{PU_{\lambda_\varepsilon,\xi_\varepsilon}^{2_{\mu}^*-1}(y)\phi(y)
\Bigg[(PU_{\lambda_\varepsilon,\xi_\varepsilon}+\phi)^{2_{\mu}^*-1}-PU_{\lambda_\varepsilon,\xi_\varepsilon}^{2_{\mu}^*-1}\Bigg](x)\varphi(x)}{|x-y|^{\mu}}dxdy
\nonumber\\
=&O\left(\|\phi\|^{\min\{2,2_{\mu}^*-1\}}\right)\|\varphi\|.
\end{align*}
\end{proof}

Now based on the above lemmas, we can prove Proposition \ref{phi}. The proof is as follows.

\noindent{\it Proof of Proposition \ref{phi}:}
It follows from Lemma \ref{Linver} that solving problem (\ref{linearprorl}) is equivalent to find a solution of
\begin{align*}
T(\phi)=\left(Q_{\varepsilon} L_{\varepsilon}\right)^{-1}(N_{\varepsilon}(\phi)+l_\varepsilon).
\end{align*}
Let
\begin{align*}
B:=\{\phi: \phi\in K^\perp, \|\phi\|\leq Cr(\varepsilon)\},
\end{align*}
where
\begin{align*}
r(\varepsilon)
=\begin{cases}
O\left(\varepsilon^{\frac{N-\mu+2}{4}}\right), \ &\text{ if } N>6-\mu,\\
O\left(\varepsilon^{\frac{N-2}{2}}\left|\ln \varepsilon\right|\right), \ &\text{ if } N=6-\mu,\\
O\left(\varepsilon^{\frac{N-2}{2}}\right), \ &\text{ if } N<6-\mu.
\end{cases}
\end{align*}
Now we shall apply the contraction mapping Theorem in a ball $B$. First of all, we prove that $T$ maps $B$ to $B$. In fact, for any $\phi\in K^\perp$,
in terms of Lemma \ref{Linver} and Lemmas \ref{lvarepsilon}, \ref{R}, we have
\begin{align*}
\|T(\phi)\|\leq C\left(\|N_{\varepsilon}(\phi)\|+\|l_\varepsilon\|\right)
\leq Cr(\varepsilon).
\end{align*}

Next, we prove that $T$ is a contraction map. To this end, choose $\phi_1,\phi_2\in K^\perp$, by Lemma \ref{Linver}, we get
\begin{align*}
\|T(\phi_1)-T(\phi_2)\|\leq C\|N_{\varepsilon}(\phi_1)-N_{\varepsilon}(\phi_2)\|.
\end{align*}
Since for $\forall \varphi\in H_0^1(\Omega_\varepsilon)$, we have
\begin{align*}
&\left\langle N_{\varepsilon}(\phi_1)-N_{\varepsilon}(\phi_2),\varphi\right\rangle
\nonumber\\
=&\int_{\Omega_\varepsilon}\Bigg\{f(PU_{\lambda_\varepsilon,\xi_\varepsilon}+\phi_1)-f(PU_{\lambda_\varepsilon,\xi_\varepsilon})
-f'(PU_{\lambda_\varepsilon,\xi_\varepsilon})\phi_1\Bigg\}\varphi
\nonumber\\
&
-\int_{\Omega_\varepsilon}\Bigg\{f(PU_{\lambda_\varepsilon,\xi_\varepsilon}+\phi_2)-f(PU_{\lambda_\varepsilon,\xi_\varepsilon})
-f'(PU_{\lambda_\varepsilon,\xi_\varepsilon})\phi_2\Bigg\}\varphi
\nonumber\\
=&\int_{\Omega_\varepsilon}\Bigg\{f(PU_{\lambda_\varepsilon,\xi_\varepsilon}+\phi_1)-f(PU_{\lambda_\varepsilon,\xi_\varepsilon}+\phi_2)
-f'(PU_{\lambda_\varepsilon,\xi_\varepsilon})(\phi_1-\phi_2)\Bigg\}\varphi
\nonumber\\
=&\int_{\Omega_\varepsilon}\int_{\Omega_\varepsilon}\frac{\left(PU_{\lambda_\varepsilon,\xi_\varepsilon}+\phi_1\right)^{2_{\mu}^*}(y)
\left(PU_{\lambda_\varepsilon,\xi_\varepsilon}+\phi_1\right)^{2_{\mu}^*-1}(x)\varphi(x)}{|x-y|^{\mu}}dxdy
\nonumber\\
&-\int_{\Omega_\varepsilon}\int_{\Omega_\varepsilon}\frac{\left(PU_{\lambda_\varepsilon,\xi_\varepsilon}+\phi_2\right)^{2_{\mu}^*}(y)
\left(PU_{\lambda_\varepsilon,\xi_\varepsilon}+\phi_2\right)^{2_{\mu}^*-1}(x)\varphi(x)}{|x-y|^{\mu}}dxdy
\nonumber\\
&-2_{\mu}^*\int_{\Omega_\varepsilon}\int_{\Omega_\varepsilon}\frac{PU_{\lambda_\varepsilon,\xi_\varepsilon}^{2_{\mu}^*-1}(y)(\phi_1-\phi_2)(y)
PU_{\lambda_\varepsilon,\xi_\varepsilon}^{2_{\mu}^*-1}(x)\varphi(x)}{|x-y|^{\mu}}dxdy
\nonumber\\
&-(2_{\mu}^*-1)\int_{\Omega_\varepsilon}\int_{\Omega_\varepsilon}\frac{PU_{\lambda_\varepsilon,\xi_\varepsilon}^{2_{\mu}^*}(y)
PU_{\lambda_\varepsilon,\xi_\varepsilon}^{2_{\mu}^*-2}(x)(\phi_1-\phi_2)(x)\varphi(x)}{|x-y|^{\mu}}dxdy
\nonumber\\
=&\int_{\Omega_\varepsilon}\int_{\Omega_\varepsilon}
\left(\left(PU_{\lambda_\varepsilon,\xi_\varepsilon}+\phi_1\right)^{2_{\mu}^*}
-\left(PU_{\lambda_\varepsilon,\xi_\varepsilon}+\phi_2\right)^{2_{\mu}^*}
-2_{\mu}^*\left(PU_{\lambda_\varepsilon,\xi_\varepsilon}+\phi_2\right)^{2_{\mu}^*-1}(\phi_1-\phi_2)
\right)(y)
\nonumber\\
&\frac{
\left(PU_{\lambda_\varepsilon,\xi_\varepsilon}+\phi_1\right)^{2_{\mu}^*-1}(x)\varphi(x)}{|x-y|^{\mu}}dxdy
\nonumber\\
&-\int_{\Omega_\varepsilon}\int_{\Omega_\varepsilon}
\left(\left(PU_{\lambda_\varepsilon,\xi_\varepsilon}+\phi_2\right)^{2_{\mu}^*-1}
-\left(PU_{\lambda_\varepsilon,\xi_\varepsilon}+\phi_1\right)^{2_{\mu}^*-1}
-(2_{\mu}^*-1)\left(PU_{\lambda_\varepsilon,\xi_\varepsilon}+\phi_1\right)^{2_{\mu}^*-2}(\phi_2-\phi_1)
\right)(x)
\nonumber\\
&\frac{\left(PU_{\lambda_\varepsilon,\xi_\varepsilon}+\phi_2\right)^{2_{\mu}^*}(y)
\varphi(x)}{|x-y|^{\mu}}dxdy
\nonumber\\
&+2_{\mu}^*\int_{\Omega_\varepsilon}\int_{\Omega_\varepsilon}\frac{\left(
\left(PU_{\lambda_\varepsilon,\xi_\varepsilon}+\phi_2\right)^{2_{\mu}^*-1}
-PU_{\lambda_\varepsilon,\xi_\varepsilon}^{2_{\mu}^*-1}
\right)(y)(\phi_1-\phi_2)(y)
\left(PU_{\lambda_\varepsilon,\xi_\varepsilon}+\phi_1\right)^{2_{\mu}^*-1}(x)\varphi(x)}{|x-y|^{\mu}}dxdy
\nonumber\\
&+2_{\mu}^*\int_{\Omega_\varepsilon}\int_{\Omega_\varepsilon}\frac{PU_{\lambda_\varepsilon,\xi_\varepsilon}^{2_{\mu}^*-1}(y)
(\phi_1-\phi_2)(y)
\left(\left(PU_{\lambda_\varepsilon,\xi_\varepsilon}+\phi_1\right)^{2_{\mu}^*-1}(x)
-PU_{\lambda_\varepsilon,\xi_\varepsilon}^{2_{\mu}^*-1}(x)
\right)\varphi(x)}{|x-y|^{\mu}}dxdy
\nonumber\\
&+(2_{\mu}^*-1)\int_{\Omega_\varepsilon}\int_{\Omega_\varepsilon}\frac{\left(
\left(PU_{\lambda_\varepsilon,\xi_\varepsilon}+\phi_2\right)^{2_{\mu}^*}
-PU_{\lambda_\varepsilon,\xi_\varepsilon}^{2_{\mu}^*}
\right)(y)
\left(PU_{\lambda_\varepsilon,\xi_\varepsilon}+\phi_1\right)^{2_{\mu}^*-2}(x)
(\phi_1-\phi_2)(x)\varphi(x)}{|x-y|^{\mu}}dxdy
\nonumber\\
&+(2_{\mu}^*-1)\int_{\Omega_\varepsilon}\int_{\Omega_\varepsilon}\frac{
PU_{\lambda_\varepsilon,\xi_\varepsilon}^{2_{\mu}^*}(y)\left(
\left(PU_{\lambda_\varepsilon,\xi_\varepsilon}+\phi_1\right)^{2_{\mu}^*-2}(x)
-PU_{\lambda_\varepsilon,\xi_\varepsilon}^{2_{\mu}^*-2}(x)
\right)
(\phi_1-\phi_2)(x)\varphi(x)}{|x-y|^{\mu}}dxdy
\nonumber\\
=&o\left(\|\phi_1-\phi_2\|\|\varphi\|\right).
\end{align*}
This implies that there exists $t\in (0,1)$ such that
\begin{align*}
\|N_{\varepsilon}(\phi_1)-N_{\varepsilon}(\phi_2)\|
\leq t\|\phi_1-\phi_2\|.
\end{align*}
Therefore,
\begin{align*}
\|T(\phi_1)-T(\phi_2)\|\leq t\|\phi_1-\phi_2\|.
\end{align*}
By the contraction mapping Theorem, we conclude that it has a unique fixed point $\phi$ in $B$.
\qed

\section{Expansion of the energy functional}\label{energy expansion}

This part is devoted to the proof of Proposition \ref{criticalpoint}.

\noindent {\it Proof of Proposition \ref{criticalpoint} (i)}:
By (\ref{linearprorl}), we get that if $(\bar{\tau},\bar{\lambda})$ is a critical point of $\tilde{I}_\varepsilon$,
\begin{align*}
\nabla I_\varepsilon(\bar{\tau},\bar{\lambda})
=I'_\varepsilon(PU_{\lambda_\varepsilon, \xi_\varepsilon}+\phi)[\nabla PU_{\lambda_\varepsilon, \xi_\varepsilon}+\nabla\phi]
=\sum_{l=0}^Nc^{l}\left\langle PZ^l_{\lambda_\varepsilon, \xi_\varepsilon},\nabla PU_{\lambda_\varepsilon, \xi_\varepsilon}+\nabla\phi\right\rangle
=0.
\end{align*}
Next, we prove that $c^{l}=0$ for $\forall l=0,\cdots, N$. Since
\begin{align*}
\partial_{\lambda}PU_{\lambda_\varepsilon, \xi_\varepsilon}
=\varepsilon^{-\frac{1}{2}}PZ_{\lambda_\varepsilon, \xi_\varepsilon}^0
-\frac{\varepsilon^{\frac{1}{2}}}{\lambda^2}\sum_{j=1}^NPZ_{\lambda_\varepsilon, \xi_\varepsilon}^j\tau^{j},\quad
\nabla_{\tau^r}PU_{\lambda_\varepsilon, \xi_\varepsilon}=\lambda^{-1}_{\varepsilon}PZ_{\lambda_\varepsilon, \xi_\varepsilon}^{r},
\end{align*}
we have
\begin{align*}
\sum_{l=0}^Nc^{l}\left\langle PZ^l_{\lambda_\varepsilon, \xi_\varepsilon},\partial_{\lambda}
PU_{\lambda_\varepsilon, \xi_\varepsilon}\right\rangle
=&\varepsilon^{-\frac{1}{2}}\sum_{l=0}^Nc^{l}\left\langle PZ^l_{\lambda_\varepsilon, \xi_\varepsilon}, PZ_{\lambda_\varepsilon, \xi_\varepsilon}^0\right\rangle
-\frac{\varepsilon^{\frac{1}{2}}}{\lambda^2}\sum_{l=0}^N\sum_{j=1}^Nc^{l}\tau^{j}\left\langle PZ^l_{\lambda_\varepsilon, \xi_\varepsilon}, PZ_{\lambda_\varepsilon, \xi_\varepsilon}^j\right\rangle
\nonumber\\
=&\varepsilon^{-\frac{1}{2}}\left(c^0c_0+o(1)\right)
-\frac{\varepsilon^{\frac{1}{2}}}{\lambda^2}\lambda_\varepsilon^2\left(\sum_{j=1}^Nc^{j}\tau^{j}c_j(1+o(1))\right)
\nonumber\\
=&\varepsilon^{-\frac{1}{2}}\left(c^0c_0+o(1)\right)
-\varepsilon^{-\frac{1}{2}}\left(\sum_{j=1}^Nc^{j}\tau^{j}c_j(1+o(1))\right)
\end{align*}
and
\begin{align*}
\sum_{l=0}^Nc^{l}\left\langle PZ^l_{\lambda_\varepsilon, \xi_\varepsilon},\partial_{\tau^r}PU_{\lambda_\varepsilon, \xi_\varepsilon}\right\rangle
=\lambda^{-1}_{\varepsilon}\sum_{l=0}^Nc^{l}\left\langle PZ^l_{\lambda_\varepsilon, \xi_\varepsilon},
PZ_{\lambda_\varepsilon, \xi_\varepsilon}^{r}\right\rangle
=\lambda_{\varepsilon}c^{r}(c_r+o(1)).
\end{align*}
On the other hand,
\begin{align*}
\left\langle PZ^l_{\lambda_\varepsilon, \xi_\varepsilon},\partial_s\phi\right\rangle
=-\left\langle \partial_sPZ^l_{\lambda_\varepsilon, \xi_\varepsilon},\phi\right\rangle
=O\left(\|\partial_sPZ^l_{\lambda_\varepsilon, \xi_\varepsilon}\|\|\phi\|\right)
=o\left(\|\partial_sPZ^l_{\lambda_\varepsilon, \xi_\varepsilon}\|\right).
\end{align*}
Then we use the fact that
\begin{align*}
\|\partial_sPZ_{\lambda_\varepsilon, \xi_\varepsilon}^l\|=
\begin{cases}
  \|\partial_{\lambda}PZ_{\lambda_\varepsilon, \xi_\varepsilon}^l\|=O(\varepsilon^{-\frac{1}{2}}), \\
  \|\partial_{\tau^j}PZ_{\lambda_\varepsilon, \xi_\varepsilon}^l\|=O(\lambda_{\varepsilon}),
  \end{cases}
\end{align*}
we can deduce that $c^l=0$ for all $l=0,\cdots,N$.
\qed

The following estimate is crucial to find critical points of $I_\varepsilon$.

\begin{lemma}\label{IVphic1}
For any $\eta>0$ small but fixed, there exists $\varepsilon_0>0$ such that for any $\varepsilon\in (0,\varepsilon_0)$, we have
\begin{align*}
I_\varepsilon(PU_{\lambda_\varepsilon, \xi_\varepsilon}+\phi)
=I_\varepsilon(PU_{\lambda_\varepsilon, \xi_\varepsilon})
+o\left(\varepsilon^{\frac{N-2}{2}}\right),
\end{align*}
$C^1$-uniformly with respect to $\lambda$ and $\tau$, satisfying (\ref{lambdajtauj}).
\end{lemma}
\begin{proof}
Since
\begin{align*}
&I_\varepsilon(PU_{\lambda_\varepsilon, \xi_\varepsilon}+\phi)-I_\varepsilon(PU_{\lambda_\varepsilon, \xi_\varepsilon})
\nonumber\\
=&\frac{1}{2}\|\phi\|^2-\Bigg[\int_{\Omega_\varepsilon}f(PU_{\lambda_\varepsilon, \xi_\varepsilon})\phi-\int_{\Omega_\varepsilon}f(U_{\lambda_\varepsilon, \xi_\varepsilon})\phi\Bigg]
\nonumber\\
&-\int_{\Omega_\varepsilon}\Bigg(F(PU_{\lambda_\varepsilon, \xi_\varepsilon}+\phi)-F(PU_{\lambda_\varepsilon, \xi_\varepsilon})-f(PU_{\lambda_\varepsilon, \xi_\varepsilon})\phi
\Bigg)
\nonumber\\
=&\frac{1}{2}\|\phi\|^2-\left\langle l_\varepsilon,\phi\right\rangle
-\int_{\Omega_\varepsilon}\Bigg(F(PU_{\lambda_\varepsilon, \xi_\varepsilon}+\phi)
-F(PU_{\lambda_\varepsilon, \xi_\varepsilon})-f(PU_{\lambda_\varepsilon, \xi_\varepsilon})\phi\Bigg),
\end{align*}
where $F(u):=\frac{1}{2\cdot 2_{\mu}^*}\left(\int_{\Omega_\varepsilon}
\frac{u^{2_\mu^*}(y)}{|x-y|^{\mu}}dy\right)u^{2_\mu^*}(x)$. Then
\begin{align*}
&\int_{\Omega_\varepsilon}\Bigg(F(PU_{\lambda_\varepsilon, \xi_\varepsilon}+\phi)
-F(PU_{\lambda_\varepsilon, \xi_\varepsilon})-f(PU_{\lambda_\varepsilon, \xi_\varepsilon})\phi\Bigg)
\nonumber\\
=&\frac{1}{2\cdot2_\mu^*}\int_{\Omega_\varepsilon}\int_{\Omega_\varepsilon}\frac{\Bigg((PU_{\lambda_\varepsilon, \xi_\varepsilon}+\phi)^{2_\mu^*}
-PU_{\lambda_\varepsilon, \xi_\varepsilon}^{2_\mu^*}
-2_\mu^*PU_{\lambda_\varepsilon, \xi_\varepsilon}^{2_\mu^*-1}(y)\phi
\Bigg)(y)
(PU_{\lambda_\varepsilon, \xi_\varepsilon}+\phi)^{2_\mu^*}(x)}{|x-y|^{\mu}}dxdy
\nonumber\\
&+\frac{1}{2\cdot2_\mu^*}\int_{\Omega_\varepsilon}\int_{\Omega_\varepsilon}\frac{\Bigg((PU_{\lambda_\varepsilon, \xi_\varepsilon}+\phi)^{2_\mu^*}(x)
-PU_{\lambda_\varepsilon, \xi_\varepsilon}^{2_\mu^*}(x)
-2_\mu^*PU_{\lambda_\varepsilon, \xi_\varepsilon}^{2_\mu^*-1}(x)\phi(x)
\Bigg)
PU_{\lambda_\varepsilon, \xi_\varepsilon}^{2_\mu^*}(y)
}{|x-y|^{\mu}}dxdy
\nonumber\\
&+\frac{1}{2}\int_{\Omega_\varepsilon}\int_{\Omega_\varepsilon}\frac{\Bigg((PU_{\lambda_\varepsilon, \xi_\varepsilon}+\phi)^{2_\mu^*}(x)
-PU_{\lambda_\varepsilon, \xi_\varepsilon}^{2_\mu^*}(x)
\Bigg)
PU_{\lambda_\varepsilon, \xi_\varepsilon}^{2_\mu^*-1}(y)\phi(y)}
{|x-y|^{\mu}}dxdy
\nonumber\\
=&O\left(\|\phi\|^2\right).
\end{align*}
Thus, when $0<\mu<4$, in view of Proposition \ref{phi}, we get
\begin{align*}
I_\varepsilon(PU_{\lambda_\varepsilon, \xi_\varepsilon}+\phi)-I_\varepsilon(PU_{\lambda_\varepsilon, \xi_\varepsilon})
=O\left(\|\phi\|^2\right)
=o\left(\varepsilon^{\frac{N-2}{2}}\right).
\end{align*}
Next, we prove that
\begin{align*}
\nabla I_\varepsilon(PU_{\lambda_\varepsilon, \xi_\varepsilon}+\phi)
=\nabla I_\varepsilon(PU_{\lambda_\varepsilon, \xi_\varepsilon})+o\left(\varepsilon^{\frac{N-2}{2}}\right).
\end{align*}
It holds
\begin{align}\label{gradient}
\nabla I_\varepsilon(PU_{\lambda_\varepsilon, \xi_\varepsilon}+\phi)-\nabla I_\varepsilon(PU_{\lambda_\varepsilon, \xi_\varepsilon})
=&\left(I'_\varepsilon(PU_{\lambda_\varepsilon, \xi_\varepsilon}+\phi)
-I'_\varepsilon(PU_{\lambda_\varepsilon, \xi_\varepsilon})\right)[\nabla PU_{\lambda_\varepsilon, \xi_\varepsilon}]
\nonumber\\
&+I'_\varepsilon(PU_{\lambda_\varepsilon, \xi_\varepsilon}+\phi)[\nabla \phi].
\end{align}
Let $\partial_s$ denote $\partial_{\lambda}$ and $\partial_{\tau^i}$ for $i=1,\cdots,N$. We know that $\partial_{s} PU_{\lambda_\varepsilon, \xi_\varepsilon}$ is a linear combination of $\frac{1}{\sqrt{\varepsilon}}PZ_{\lambda_\varepsilon, \xi_\varepsilon}^h$ for $h=0$ and $\sqrt{\varepsilon}PZ_{\lambda_\varepsilon, \xi_\varepsilon}^j$ for $j=1,\cdots,N$ with coefficients uniformly bounded as $\varepsilon\rightarrow 0$
for any $\lambda,\tau$ satisfying (\ref{lambdajtauj}). In order to prove (\ref{gradient}), it is enough to estimate $\left(I'_\varepsilon(PU_{\lambda_\varepsilon, \xi_\varepsilon}+\phi)
-I'_\varepsilon(PU_{\lambda_\varepsilon, \xi_\varepsilon})\right)[\frac{1}{\sqrt{\varepsilon}}PZ_{\lambda_\varepsilon, \xi_\varepsilon}^h]$ and $I'_\varepsilon(PU_{\lambda_\varepsilon, \xi_\varepsilon}+\phi)[\nabla \phi]$ because $\left(I'_\varepsilon(PU_{\lambda_\varepsilon, \xi_\varepsilon}+\phi)
-I'_\varepsilon(PU_{\lambda_\varepsilon, \xi_\varepsilon})\right)[\sqrt{\varepsilon}PZ_{\lambda_\varepsilon, \xi_\varepsilon}^j]$ can be estimated in the same way. Thus,
\begin{align*}
&\left(I'_\varepsilon(PU_{\lambda_\varepsilon, \xi_\varepsilon}+\phi)
-I'_\varepsilon(PU_{\lambda_\varepsilon, \xi_\varepsilon})\right)[\frac{1}{\sqrt{\varepsilon}}PZ_{\lambda_\varepsilon, \xi_\varepsilon}^h]
\nonumber\\
=&-\int_{\Omega_\varepsilon}\Bigg[f(PU_{\lambda_\varepsilon, \xi_\varepsilon}+\phi)-f(PU_{\lambda_\varepsilon, \xi_\varepsilon})
\Bigg]\frac{1}{\sqrt{\varepsilon}}PZ_{\lambda_\varepsilon, \xi_\varepsilon}^h
\nonumber\\
=&
-\int_{\Omega_\varepsilon}\Bigg[f(PU_{\lambda_\varepsilon, \xi_\varepsilon}+\phi)-f(PU_{\lambda_\varepsilon, \xi_\varepsilon})-f'(PU_{\lambda_\varepsilon, \xi_\varepsilon})\phi
\Bigg]\frac{1}{\sqrt{\varepsilon}}PZ^h
\nonumber\\
&-\int_{\Omega_\varepsilon}f'(PU_{\lambda_\varepsilon, \xi_\varepsilon})\phi\frac{1}{\sqrt{\varepsilon}}\Bigg[PZ_{\lambda_\varepsilon, \xi_\varepsilon}^h-Z_{\lambda_\varepsilon, \xi_\varepsilon}^h
\Bigg]
-\int_{\Omega_\varepsilon}\Bigg[f'(PU_{\lambda_\varepsilon, \xi_\varepsilon})-f'(U_{\lambda_\varepsilon, \xi_\varepsilon})
\Bigg]\phi\frac{1}{\sqrt{\varepsilon}}Z_{\lambda_\varepsilon, \xi_\varepsilon}^h.
\end{align*}
First of all, we have
\begin{align*}
&\int_{\Omega_\varepsilon}\Bigg[f'(PU_{\lambda_\varepsilon, \xi_\varepsilon})-f'(U_{\lambda_\varepsilon, \xi_\varepsilon})
\Bigg]\phi\frac{1}{\sqrt{\varepsilon}}Z_{\lambda_\varepsilon, \xi_\varepsilon}^h
\nonumber\\
=&2_{\mu}^*\Bigg\{\int_{\Omega_\varepsilon}\int_{\Omega_\varepsilon}
\frac{\left(PU_{\lambda_\varepsilon, \xi_\varepsilon}^{2_{\mu}^*-1}
-U_{\lambda_\varepsilon, \xi_\varepsilon}^{2_{\mu}^*-1}\right)
\phi(y)PU_{\lambda_\varepsilon, \xi_\varepsilon}^{2_{\mu}^*-1} \frac{1}{\sqrt{\varepsilon}}Z_{\lambda_\varepsilon, \xi_\varepsilon}^h(x)
}{|x-y|^{\mu}}dxdy
\Bigg\}
\nonumber\\
&+2_{\mu}^*\Bigg\{\int_{\Omega_\varepsilon}\int_{\Omega_\varepsilon}
\frac{U_{\lambda_\varepsilon, \xi_\varepsilon}^{2_{\mu}^*-1}\phi(y)\left(PU_{\lambda_\varepsilon, \xi_\varepsilon}^{2_{\mu}^*-1}
-U_{\lambda_\varepsilon, \xi_\varepsilon}^{2_{\mu}^*-1}
\right) \frac{1}{\sqrt{\varepsilon}}Z_{\lambda_\varepsilon, \xi_\varepsilon}^h(x)
}{|x-y|^{\mu}}dxdy
\Bigg\}
\nonumber\\
&+(2_{\mu}^*-1)\Bigg\{\int_{\Omega_\varepsilon}\int_{\Omega_\varepsilon}
\frac{\left(PU_{\lambda_\varepsilon, \xi_\varepsilon}^{2_{\mu}^*}
-U_{\lambda_\varepsilon, \xi_\varepsilon}^{2_{\mu}^*}
\right)(y)PU_{\lambda_\varepsilon, \xi_\varepsilon}^{2_{\mu}^*-2}\phi\frac{1}{\sqrt{\varepsilon}}Z_{\lambda_\varepsilon, \xi_\varepsilon}^h(x)
}{|x-y|^{\mu}}dxdy
\Bigg\}
\nonumber\\
&+(2_{\mu}^*-1)\Bigg\{\int_{\Omega_\varepsilon}\int_{\Omega_\varepsilon}
\frac{U_{\lambda_\varepsilon, \xi_\varepsilon}^{2_{\mu}^*}(y)\left(PU_{\lambda_\varepsilon, \xi_\varepsilon}^{2_{\mu}^*-2}
-U_{\lambda_\varepsilon, \xi_\varepsilon}^{2_{\mu}^*-2}
\right)\phi\frac{1}{\sqrt{\varepsilon}}Z_{\lambda_\varepsilon, \xi_\varepsilon}^h(x)
}{|x-y|^{\mu}}dxdy
\Bigg\}
\nonumber\\
=&O\left(\left(\int_{\Omega_\varepsilon}\left|PU_{\lambda_\varepsilon, \xi_\varepsilon}^{2_{\mu}^*-1}
-U_{\lambda_\varepsilon, \xi_\varepsilon}^{2_{\mu}^*-1}
\right|^{\frac{2N}{N-\mu+2}}\right)^{\frac{N-\mu+2}{2N}}
\right)\|\phi\|
+O\left(\left(\int_{\Omega_\varepsilon}\left|PU_{\lambda_\varepsilon, \xi_\varepsilon}^{2_{\mu}^*}
-U_{\lambda_\varepsilon, \xi_\varepsilon}^{2_{\mu}^*}
\right|^{\frac{2N}{2N-\mu}}\right)^{\frac{2N-\mu}{2N}}
\right)\|\phi\|
\nonumber\\
&+O\left(\left(\int_{\Omega_\varepsilon}\left|PU_{\lambda_\varepsilon, \xi_\varepsilon}^{2_{\mu}^*-2}
-U_{\lambda_\varepsilon, \xi_\varepsilon}^{2_{\mu}^*-2}
\right|^{\frac{2N}{4-\mu}}\right)^{\frac{4-\mu}{2N}}
\right)\|\phi\|
\nonumber\\
=&o\left(\varepsilon^{\frac{N-2}{2}}\right).
\end{align*}
Moreover, by Lemmas \ref{PU} and \ref{lvarepsilon},
\begin{align*}
\int_{\Omega_\varepsilon}f'(PU_{\lambda_\varepsilon, \xi_\varepsilon})\phi\frac{1}{\sqrt{\varepsilon}}
\Bigg[PZ_{\lambda_\varepsilon, \xi_\varepsilon}^h-Z_{\lambda_\varepsilon, \xi_\varepsilon}^h
\Bigg]
=O\left(\frac{1}{\sqrt{\varepsilon}}\left(\int_{\Omega_\varepsilon}\left|PZ_{\lambda_\varepsilon, \xi_\varepsilon}^h-Z_{\lambda_\varepsilon, \xi_\varepsilon}^h
\right|^{\frac{2N}{N-2}}\right)^{\frac{N-2}{2N}}
\right)\|\phi\|
=o\left(\varepsilon^{\frac{N-2}{2}}\right),
\end{align*}
and because of $PU_{\lambda_\varepsilon, \xi_\varepsilon}\leq U_{\lambda_\varepsilon, \xi_\varepsilon}\leq C\varepsilon^{\frac{N-2}{4}}$ in $\Omega_\varepsilon\setminus B(0,\rho)$, Hardy-Littlewood-Sobolev inequality, H\"{o}lder inequality and Lemma \ref{lvarepsilon}, we get
\begin{align*}
&\int_{\Omega_\varepsilon\setminus B(0,\rho)}
\Bigg[f(PU_{\lambda_\varepsilon, \xi_\varepsilon}+\phi)-f(PU_{\lambda_\varepsilon, \xi_\varepsilon})-f'(PU_{\lambda_\varepsilon, \xi_\varepsilon})\phi
\Bigg]\frac{1}{\sqrt{\varepsilon}}PZ_{\lambda_\varepsilon, \xi_\varepsilon}^h
\nonumber\\
=&O\left(\left(\int_{\Omega}\left(PU_{\lambda_\varepsilon, \xi_\varepsilon}+\phi
\right)^{2_{\mu}^*\frac{2N}{2N-\mu}}\right)^{\frac{2N-\mu}{2N}}
\left(\int_{\Omega_\varepsilon\setminus B(0,\rho)}\left(
PU_{\lambda_\varepsilon, \xi_\varepsilon}+\phi\right)^{(2_{\mu}^*-1)\frac{2N}{2N-\mu}}
\left(\frac{1}{\sqrt{\varepsilon}}PZ_{\lambda_\varepsilon, \xi_\varepsilon}^h\right)^{\frac{2N}{2N-\mu}}
\right)^{\frac{2N-\mu}{2N}}
\right)
\nonumber\\
&+O\left(\left(\int_{\Omega}PU_{\lambda_\varepsilon, \xi_\varepsilon}^{2_{\mu}^*\frac{2N}{2N-\mu}}\right)^{\frac{2N-\mu}{2N}}
\left(\int_{\Omega_\varepsilon\setminus B(0,\rho)}
PU_{\lambda_\varepsilon, \xi_\varepsilon}^{(2_{\mu}^*-1)\frac{2N}{2N-\mu}}
\left(\frac{1}{\sqrt{\varepsilon}}PZ_{\lambda_\varepsilon, \xi_\varepsilon}^h\right)^{\frac{2N}{2N-\mu}}
\right)^{\frac{2N-\mu}{2N}}
\right)
\nonumber\\
&+O\left(\left(\int_{\Omega}\left(PU_{\lambda_\varepsilon, \xi_\varepsilon}^{2_{\mu}^*-1}\phi
\right)^{\frac{2N}{2N-\mu}}\right)^{\frac{2N-\mu}{2N}}
\left(\int_{\Omega_\varepsilon\setminus B(0,\rho)}\left(
\frac{1}{\sqrt{\varepsilon}}PZ_{\lambda_\varepsilon, \xi_\varepsilon}^hPU_{\lambda_\varepsilon, \xi_\varepsilon}^{2_{\mu}^*-1}
\right)^{\frac{2N}{2N-\mu}}\right)^{\frac{2N-\mu}{2N}}
\right)
\nonumber\\
&+O\left(\left(\int_{\Omega}\left(PU_{\lambda_\varepsilon, \xi_\varepsilon}^{2_{\mu}^*}
\right)^{\frac{2N}{2N-\mu}}\right)^{\frac{2N-\mu}{2N}}
\left(\int_{\Omega_\varepsilon\setminus B(0,\rho)}\left(
\frac{1}{\sqrt{\varepsilon}}PZ_{\lambda_\varepsilon, \xi_\varepsilon}^hPU_{\lambda_\varepsilon, \xi_\varepsilon}^{2_{\mu}^*-2}\phi
\right)^{\frac{2N}{2N-\mu}}\right)^{\frac{2N-\mu}{2N}}
\right)
=o\left(\varepsilon^{\frac{N-2}{2}}\right).
\end{align*}
Thus,
\begin{align*}
&\int_{\Omega_\varepsilon}\Bigg[f(PU_{\lambda_\varepsilon, \xi_\varepsilon}+\phi)-f(PU_{\lambda_\varepsilon, \xi_\varepsilon})-f'(PU_{\lambda_\varepsilon, \xi_\varepsilon})\phi
\Bigg]\frac{1}{\sqrt{\varepsilon}}PZ_{\lambda_\varepsilon, \xi_\varepsilon}^h
\nonumber\\
=&\Bigg\{\int_{\Omega_\varepsilon\setminus B(0,\rho)}+\int_{B(0,\rho)\setminus B(0,\varepsilon)}\Bigg\}
\Bigg[f(PU_{\lambda_\varepsilon, \xi_\varepsilon}+\phi)-f(PU_{\lambda_\varepsilon, \xi_\varepsilon})-f'(PU_{\lambda_\varepsilon, \xi_\varepsilon})\phi
\Bigg]\frac{1}{\sqrt{\varepsilon}}PZ_{\lambda_\varepsilon, \xi_\varepsilon}^h
\nonumber\\
=&\int_{B(0,\rho)\setminus B(0,\varepsilon)}\Bigg[f(PU_{\lambda_\varepsilon, \xi_\varepsilon}+\phi)-f(PU_{\lambda_\varepsilon, \xi_\varepsilon})-f'(U_{\lambda_\varepsilon, \xi_\varepsilon})\phi
\Bigg]\frac{1}{\sqrt{\varepsilon}}PZ_{\lambda_\varepsilon, \xi_\varepsilon}^h
\nonumber\\
&+\int_{B(0,\rho)\setminus B(0,\varepsilon)}\Bigg[f'(U_{\lambda_\varepsilon, \xi_\varepsilon})-f'(PU_{\lambda_\varepsilon, \xi_\varepsilon})
\Bigg]\phi\frac{1}{\sqrt{\varepsilon}}PZ_{\lambda_\varepsilon, \xi_\varepsilon}^h
+o\left(\varepsilon^{\frac{N-2}{2}}\right),
\end{align*}
where for $t\in [0,1]$,
\begin{align*}
&\int_{B(0,\rho)\setminus B(0,\varepsilon)}\Bigg[f(PU_{\lambda_\varepsilon, \xi_\varepsilon}+\phi)-f(PU_{\lambda_\varepsilon, \xi_\varepsilon})-f'(U_{\lambda_\varepsilon, \xi_\varepsilon})\phi
\Bigg]\frac{1}{\sqrt{\varepsilon}}PZ_{\lambda_\varepsilon, \xi_\varepsilon}^h
\nonumber\\
=&\int_{B(0,\rho)\setminus B(0,\varepsilon)}\Bigg[f'(PU_{\lambda_\varepsilon, \xi_\varepsilon}+t\phi)-f'(U_{\lambda_\varepsilon, \xi_\varepsilon})
\Bigg]\phi\frac{1}{\sqrt{\varepsilon}}PZ_{\lambda_\varepsilon, \xi_\varepsilon}^h
\nonumber\\
=&\int_{\Omega_\varepsilon}\int_{B(0,\rho)\setminus B(0,\varepsilon)}\frac{2_{\mu}^*(PU_{\lambda_\varepsilon, \xi_\varepsilon}+t\phi)^{2_{\mu}^*-1}\phi(y)
(PU_{\lambda_\varepsilon, \xi_\varepsilon}+t\phi)^{2_{\mu}^*-1}\frac{1}{\sqrt{\varepsilon}}PZ_{\lambda_\varepsilon, \xi_\varepsilon}^h
}{|x-y|^{\mu}}dxdy
\nonumber\\
&+\int_{\Omega_\varepsilon}\int_{B(0,\rho)\setminus B(0,\varepsilon)}\frac{(2_{\mu}^*-1)(PU_{\lambda_\varepsilon, \xi_\varepsilon}+t\phi)^{2_{\mu}^*}(y)
(PU_{\lambda_\varepsilon, \xi_\varepsilon}+t\phi)^{2_{\mu}^*-2}\phi\frac{1}{\sqrt{\varepsilon}}PZ_{\lambda_\varepsilon, \xi_\varepsilon}^h(x)
}{|x-y|^{\mu}}dxdy
\nonumber\\
&-\Bigg(\int_{\Omega_\varepsilon}\int_{B(0,\rho)\setminus B(0,\varepsilon)}\frac{2_{\mu}^*U_{\lambda_\varepsilon, \xi_\varepsilon}^{2_{\mu}^*-1}\phi(y)
U_{\lambda_\varepsilon, \xi_\varepsilon}^{2_{\mu}^*-1}\frac{1}{\sqrt{\varepsilon}}PZ_{\lambda_\varepsilon, \xi_\varepsilon}^h
}{|x-y|^{\mu}}dxdy
\nonumber\\
&+\int_{\Omega_\varepsilon}\int_{B(0,\rho)\setminus B(0,\varepsilon)}\frac{(2_{\mu}^*-1)U_{\lambda_\varepsilon, \xi_\varepsilon}^{2_{\mu}^*}(y)
U_{\lambda_\varepsilon, \xi_\varepsilon}^{2_{\mu}^*-2}\phi\frac{1}{\sqrt{\varepsilon}}PZ_{\lambda_\varepsilon, \xi_\varepsilon}^h(x)
}{|x-y|^{\mu}}dxdy
\Bigg)
\nonumber\\
=&2_{\mu}^*\int_{\Omega_\varepsilon}\int_{B(0,\rho)\setminus B(0,\varepsilon)}\frac{\Bigg((PU_{\lambda_\varepsilon, \xi_\varepsilon}+t\phi)^{2_{\mu}^*-1}-U_{\lambda_\varepsilon, \xi_\varepsilon}^{2_{\mu}^*-1}\Bigg)\phi(y)
(PU_{\lambda_\varepsilon, \xi_\varepsilon}+t\phi)^{2_{\mu}^*-1}\frac{1}{\sqrt{\varepsilon}}PZ_{\lambda_\varepsilon, \xi_\varepsilon}^h(x)
}{|x-y|^{\mu}}dxdy
\nonumber\\
&+2_{\mu}^*\int_{\Omega_\varepsilon}\int_{B(0,\rho)\setminus B(0,\varepsilon)}
\frac{U_{\lambda_\varepsilon, \xi_\varepsilon}^{2_{\mu}^*-1}\phi(y)
\Bigg((PU_{\lambda_\varepsilon, \xi_\varepsilon}+t\phi)^{2_{\mu}^*-1}-U_{\lambda_\varepsilon, \xi_\varepsilon}^{2_{\mu}^*-1}
\Bigg)
\frac{1}{\sqrt{\varepsilon}}PZ_{\lambda_\varepsilon, \xi_\varepsilon}^h(x)
}{|x-y|^{\mu}}dxdy
\nonumber\\
&+(2_{\mu}^*-1)\int_{\Omega_\varepsilon}\int_{B(0,\rho)\setminus B(0,\varepsilon)}\frac{\Bigg((PU_{\lambda_\varepsilon, \xi_\varepsilon}+t\phi)^{2_{\mu}^*}
-U_{\lambda_\varepsilon, \xi_\varepsilon}^{2_{\mu}^*}\Bigg)(y)
(PU_{\lambda_\varepsilon, \xi_\varepsilon}+t\phi)^{2_{\mu}^*-2}\phi\frac{1}{\sqrt{\varepsilon}}PZ_{\lambda_\varepsilon, \xi_\varepsilon}^h(x)
}{|x-y|^{\mu}}dxdy
\nonumber\\
&+(2_{\mu}^*-1)\int_{\Omega_\varepsilon}\int_{B(0,\rho)\setminus B(0,\varepsilon)}
\frac{U_{\lambda_\varepsilon, \xi_\varepsilon}^{2_{\mu}^*}(y)\Bigg((PU_{\lambda_\varepsilon, \xi_\varepsilon}+t\phi)^{2_{\mu}^*-2}-U_{\lambda_\varepsilon, \xi_\varepsilon}^{2_{\mu}^*-2}\Bigg)
\phi\frac{1}{\sqrt{\varepsilon}}PZ_{\lambda_\varepsilon, \xi_\varepsilon}^h(x)
}{|x-y|^{\mu}}dxdy
\nonumber\\
=&o\left(\varepsilon^{\frac{N-2}{2}}\right).
\end{align*}
Hence, we get
\begin{align*}
\left(I'_\varepsilon(PU_{\lambda_\varepsilon, \xi_\varepsilon}+\phi)
-I'_\varepsilon(PU_{\lambda_\varepsilon, \xi_\varepsilon})\right)[\frac{1}{\sqrt{\varepsilon}}PZ_{\lambda_\varepsilon, \xi_\varepsilon}^h]
=o\left(\varepsilon^{\frac{N-2}{2}}\right).
\end{align*}
As for $I'_\varepsilon(PU_{\lambda_\varepsilon, \xi_\varepsilon})[\nabla \phi]$, we have
\begin{align*}
I'_\varepsilon(PU_{\lambda_\varepsilon, \xi_\varepsilon}+\phi)[\nabla \phi]
=\sum_{l=0}^Nc^l\left\langle PZ^l, \nabla \phi\right\rangle,
\end{align*}
and it follows from Proposition \ref{criticalpoint} (i) that $|c^l|\leq C\lambda^{-1}_{\varepsilon}\|\phi\|_{\frac{2N}{N-2}}$.
Then
\begin{align*}
\left|\left\langle PZ_{\lambda_\varepsilon, \xi_\varepsilon}^l, \nabla \phi\right\rangle\right|
\leq C\|Z_{\lambda_\varepsilon, \xi_\varepsilon}^l\|_{\frac{2N}{N-2}}\|\phi\|_{\frac{2N}{N-2}}
\leq C\lambda_{\varepsilon}\|\phi\|_{\frac{2N}{N-2}}.
\end{align*}
Furthermore,
\begin{align*}
\left|I'_\varepsilon(PU_{\lambda_\varepsilon, \xi_\varepsilon})[\nabla \phi]\right|
=O\left(\|\phi\|^2\right)
=o\left(\varepsilon^{\frac{N-2}{2}}\right).
\end{align*}
Above all, we conclude the proof.
\end{proof}

\begin{lemma}\label{IPU}
For $N\geq 5$, we have
\begin{align*}
I_\varepsilon(PU_{\lambda_\varepsilon, \xi_\varepsilon})
=&\left(1-\frac{1}{2_\mu^*}\right)\frac{N(N-2)}{2\mathcal{A}_{H,L}}A_N
+\frac{N(N-2)^2\omega_{N}B_N}{2\mathcal{A}_{H,L}}\frac{H(0,0)}{\lambda_\varepsilon^{N-2}}\left(1+o(1)\right)
\nonumber\\
&+\frac{N(N-2)}{2\mathcal{A}_{H,L}(1+|\tau|^2)^{\frac{N-2}{2}}}
M(\tau)\varepsilon^{N-2}\lambda_\varepsilon^{N-2}\left(1+o(1)\right),
\end{align*}
$C^1$-uniformly with respect to $\lambda$ and $\tau$ satisfying (\ref{lambdajtauj}).
Here
\begin{align*}
A_N=\int_{\mathbb{R}^N}U_{1, 0}^{2^*}(x)dx, \quad B_N=\int_{\mathbb{R}^N}U_{1, 0}^{2^*-1}(x)dx, \quad
M(\tau)=\int_{\mathbb{R}^N}\frac{1}{|z|^{N-2}(1+|z-\tau|^2)^{\frac{N+2}{2}}}dz.
\end{align*}
\end{lemma}
\begin{proof}
By the definition of $I_\varepsilon(u)$ and $PU_{\lambda_\varepsilon, \xi_\varepsilon}$, we get
\begin{align}\label{IPU1}
I_\varepsilon(PU_{\lambda_\varepsilon, \xi_\varepsilon})
=&\frac{1}{2}\int_{\Omega_\varepsilon}|\nabla PU_{\lambda_\varepsilon, \xi_\varepsilon}|^2
-\frac{1}{2\cdot2_\mu^*}\int_{\Omega_\varepsilon}\int_{\Omega_\varepsilon}\frac{PU_{\lambda_\varepsilon, \xi_\varepsilon}^{2_\mu^*}(y)PU_{\lambda_\varepsilon, \xi_\varepsilon}^{2_\mu^*}(x)}{|x-y|^{\mu}}dxdy
\nonumber\\
=&\frac{1}{2}\int_{\Omega_\varepsilon}\int_{\Omega_\varepsilon}\frac{U_{\lambda_\varepsilon, \xi_\varepsilon}^{2_\mu^*}(y)U_{\lambda_\varepsilon, \xi_\varepsilon}^{2_\mu^*-1}(x)PU_{\lambda_\varepsilon, \xi_\varepsilon}(x)}{|x-y|^{\mu}}dxdy
\nonumber\\
&-\frac{1}{2\cdot2_\mu^*}\int_{\Omega_\varepsilon}\int_{\Omega_\varepsilon}
\frac{PU_{\lambda_\varepsilon, \xi_\varepsilon}^{2_\mu^*}(y)PU_{\lambda_\varepsilon, \xi_\varepsilon}^{2_\mu^*}(x)}{|x-y|^{\mu}}dxdy
\nonumber\\
=&\frac{1}{2}\int_{\Omega_\varepsilon}\int_{\Omega_\varepsilon}\frac{U_{\lambda_\varepsilon, \xi_\varepsilon}^{2_\mu^*}(y)U_{\lambda_\varepsilon, \xi_\varepsilon}^{2_\mu^*}(x)}{|x-y|^{\mu}}dxdy
\nonumber\\
&+\frac{1}{2}\int_{\Omega_\varepsilon}\int_{\Omega_\varepsilon}\frac{U_{\lambda_\varepsilon, \xi_\varepsilon}^{2_\mu^*}(y)U_{\lambda_\varepsilon, \xi_\varepsilon}^{2_\mu^*-1}(x)(PU_{\lambda_\varepsilon, \xi_\varepsilon}(x)-U_{\lambda_\varepsilon, \xi_\varepsilon}(x))}{|x-y|^{\mu}}dxdy
\nonumber\\
&-\frac{1}{2\cdot2_\mu^*}\int_{\Omega_\varepsilon}\int_{\Omega_\varepsilon}
\frac{PU_{\lambda_\varepsilon, \xi_\varepsilon}^{2_\mu^*}(y)PU_{\lambda_\varepsilon, \xi_\varepsilon}^{2_\mu^*}(x)}{|x-y|^{\mu}}dxdy
\nonumber\\
=&\frac{1}{2}A+\frac{1}{2}B-\frac{1}{2\cdot2_\mu^*}C.
\end{align}
We first estimate $A$. Since
\begin{align}\label{A}
A=&\int_{\mathbb{R}^N}\int_{\mathbb{R}^N}\frac{U_{\lambda_\varepsilon, \xi_\varepsilon}^{2_\mu^*}(y)U_{\lambda_\varepsilon, \xi_\varepsilon}^{2_\mu^*}(x)}{|x-y|^{\mu}}dxdy
-2\int_{\mathbb{R}^N\setminus \Omega_\varepsilon}\int_{\Omega_\varepsilon}\frac{U_{\lambda_\varepsilon, \xi_\varepsilon}^{2_\mu^*}(y)U_{\lambda_\varepsilon, \xi_\varepsilon}^{2_\mu^*}(x)}{|x-y|^{\mu}}dxdy
\nonumber\\
&-\int_{\mathbb{R}^N\setminus \Omega_\varepsilon}\int_{\mathbb{R}^N\setminus \Omega_\varepsilon}\frac{U_{\lambda_\varepsilon, \xi_\varepsilon}^{2_\mu^*}(y)U_{\lambda_\varepsilon, \xi_\varepsilon}^{2_\mu^*}(x)}{|x-y|^{\mu}}dxdy
\nonumber\\
=&A_0-2A_1-A_2,
\end{align}
where
\begin{align}\label{A0}
A_0=\frac{N(N-2)}{\mathcal{A}_{H,L}}\int_{\mathbb{R}^N}U_{\lambda_\varepsilon, \xi_\varepsilon}^{2^*}(x)dx
=\frac{N(N-2)}{\mathcal{A}_{H,L}}\int_{\mathbb{R}^N}U_{1, 0}^{2^*}(x)dx
:=\frac{N(N-2)}{\mathcal{A}_{H,L}}A_N,
\end{align}
\begin{align}\label{A1}
A_1\leq
&\int_{\mathbb{R}^N\setminus \Omega_\varepsilon}\int_{\mathbb{R}^N}\frac{U_{\lambda_\varepsilon, \xi_\varepsilon}^{2_\mu^*}(y)U_{\lambda_\varepsilon, \xi_\varepsilon}^{2_\mu^*}(x)}{|x-y|^{\mu}}dxdy
\nonumber\\
=&\frac{N(N-2)}{\mathcal{A}_{H,L}}\int_{\mathbb{R}^N\setminus \Omega_\varepsilon}U_{\lambda_\varepsilon, \xi_\varepsilon}^{2^*}(x)dx
=O\left(\frac{1}{\lambda_\varepsilon^N}+(\varepsilon\lambda_\varepsilon)^N
\right)
=O\left(\varepsilon^{\frac{N}{2}}
\right).
\end{align}
Similarly,
\begin{align}\label{A2}
A_2=O\left(\frac{1}{\lambda_\varepsilon^N}+(\varepsilon\lambda_\varepsilon)^N
\right)
=O\left(\varepsilon^{\frac{N}{2}}
\right).
\end{align}
From (\ref{A}), (\ref{A0}), (\ref{A1}) and (\ref{A2}), we deduce
\begin{align}\label{ah}
A=\frac{N(N-2)}{\mathcal{A}_{H,L}}A_N
+O\left(\varepsilon^{\frac{N}{2}}
\right).
\end{align}
As to $B$, we find that
\begin{align*}
&\int_{\Omega_\varepsilon}\int_{\Omega_\varepsilon}\frac{U_{\lambda_\varepsilon, \xi_\varepsilon}^{2_\mu^*}(y)U_{\lambda_\varepsilon, \xi_\varepsilon}^{2_\mu^*-1}(x)}{|x-y|^{\mu}}dxdy
\nonumber\\
=&\Bigg\{\int_{\mathbb{R}^N}\int_{\mathbb{R}^N}
-2\int_{\mathbb{R}^N\setminus \Omega_\varepsilon}\int_{\Omega_\varepsilon}
-\int_{\mathbb{R}^N\setminus \Omega_\varepsilon}\int_{\mathbb{R}^N\setminus \Omega_\varepsilon}\Bigg\}
\frac{U_{\lambda_\varepsilon, \xi_\varepsilon}^{2_\mu^*}(y)U_{\lambda_\varepsilon, \xi_\varepsilon}^{2_\mu^*-1}(x)}{|x-y|^{\mu}}dxdy
\nonumber\\
=&\frac{N(N-2)}{\mathcal{A}_{H,L}\lambda_\varepsilon^{\frac{N-2}{2}}}\int_{\mathbb{R}^N}U_{1, 0}^{2^*-1}(x)dx-2B_1-B_2
=\frac{N(N-2)}{\mathcal{A}_{H,L}}\frac{B_N}{\lambda_\varepsilon^{\frac{N-2}{2}}}-2B_1-B_2.
\end{align*}
Moreover,
\begin{align*}
B_1\leq &\int_{\mathbb{R}^N\setminus \Omega_\varepsilon}\int_{\mathbb{R}^N}\frac{U_{\lambda_\varepsilon, \xi_\varepsilon}^{2_\mu^*}(y)U_{\lambda_\varepsilon, \xi_\varepsilon}^{2_\mu^*-1}(x)}{|x-y|^{\mu}}dxdy
\nonumber\\
=&\frac{N(N-2)}{\mathcal{A}_{H,L}}\int_{\mathbb{R}^N\setminus \Omega_\varepsilon}U_{\lambda_\varepsilon, \xi_\varepsilon}^{2^*-1}(x)dx
=O\left(\frac{1}{\lambda_\varepsilon^{\frac{N-2}{2}}}\left(
\frac{1}{\lambda_\varepsilon^2}+(\varepsilon\lambda_\varepsilon)^N\right)
\right)
=O\left(\varepsilon^{\frac{N-2}{4}}\left(\varepsilon
+\varepsilon^{\frac{N}{2}}
\right)
\right).
\end{align*}
In a same way,
\begin{align*}
B_2=O\left(\frac{1}{\lambda_\varepsilon^{\frac{N-2}{2}}}
\left(\frac{1}{\lambda_\varepsilon^2}+(\varepsilon\lambda_\varepsilon)^N\right)
\right)
=O\left(\varepsilon^{\frac{N-2}{4}}\left(\varepsilon
+\varepsilon^{\frac{N}{2}}
\right)
\right).
\end{align*}
Consequently,
\begin{align}\label{B0}
\int_{\Omega_\varepsilon}\int_{\Omega_\varepsilon}\frac{U_{\lambda_\varepsilon, \xi_\varepsilon}^{2_\mu^*}(y)U_{\lambda_\varepsilon, \xi_\varepsilon}^{2_\mu^*-1}(x)}{|x-y|^{\mu}}dxdy
=\frac{N(N-2)}{\mathcal{A}_{H,L}}\frac{B_N}{\lambda_\varepsilon^{\frac{N-2}{2}}}
+O\left(\varepsilon^{\frac{N-2}{4}}\left(\varepsilon
+\varepsilon^{\frac{N}{2}}
\right)
\right).
\end{align}
Using Lemma \ref{PU} and the Taylor expansion, we see
\begin{align}\label{B}
B=&\int_{\Omega_\varepsilon}\int_{\Omega_\varepsilon}\frac{U_{\lambda_\varepsilon, \xi_\varepsilon}^{2_\mu^*}(y)U_{\lambda_\varepsilon, \xi_\varepsilon}^{2_\mu^*-1}(x)(PU_{\lambda_\varepsilon, \xi_\varepsilon}(x)-U_{\lambda_\varepsilon, \xi_\varepsilon}(x))}{|x-y|^{\mu}}dxdy
\nonumber\\
=&-\int_{\Omega_\varepsilon}\int_{\Omega_\varepsilon}\frac{U_{\lambda_\varepsilon, \xi_\varepsilon}^{2_\mu^*}(y)U_{\lambda_\varepsilon, \xi_\varepsilon}^{2_\mu^*-1}(x)\left(
\frac{(N-2)\omega_{N}}{\lambda_\varepsilon^{\frac{N-2}{2}}}H(x,\xi_\varepsilon)
+\frac{\lambda_\varepsilon^{\frac{N-2}{2}}}{(1+|\tau|^2)^{\frac{N-2}{2}}}\frac{\varepsilon^{N-2}}{|x|^{N-2}}-R(x)
\right)}{|x-y|^{\mu}}dxdy
\nonumber\\
=&-\int_{\Omega_\varepsilon}\int_{\Omega_\varepsilon}\frac{(N-2)\omega_{N}}{\lambda_\varepsilon^{\frac{N-2}{2}}}\Bigg[
H(\xi_\varepsilon,\xi_\varepsilon)+\left\langle \nabla H(\xi_\varepsilon,\xi_\varepsilon), x-\xi_\varepsilon\right\rangle+O\left(|x-\xi_\varepsilon|^2\right)
\Bigg]
\frac{U_{\lambda_\varepsilon, \xi_\varepsilon}^{2_\mu^*}(y)U_{\lambda_\varepsilon, \xi_\varepsilon}^{2_\mu^*-1}(x)}{|x-y|^{\mu}}dxdy
\nonumber\\
&-\frac{\varepsilon^{N-2}\lambda_\varepsilon^{\frac{N-2}{2}}}{(1+|\tau|^2)^{\frac{N-2}{2}}}
\int_{\Omega_\varepsilon}\int_{\Omega_\varepsilon}\frac{1}{|x|^{N-2}}
\frac{U_{\lambda_\varepsilon, \xi_\varepsilon}^{2_\mu^*}(y)U_{\lambda_\varepsilon, \xi_\varepsilon}^{2_\mu^*-1}(x)}{|x-y|^{\mu}}dxdy
+\int_{\Omega_\varepsilon}\int_{\Omega_\varepsilon}\frac{U_{\lambda_\varepsilon, \xi_\varepsilon}^{2_\mu^*}(y)
U_{\lambda_\varepsilon, \xi_\varepsilon}^{2_\mu^*-1}(x)R(x)}{|x-y|^{\mu}}dxdy.
\end{align}
By Taylor expansion and (\ref{B0}), we have
\begin{align}\label{B1}
&-\int_{\Omega_\varepsilon}\int_{\Omega_\varepsilon}\frac{(N-2)\omega_{N}}{\lambda_\varepsilon^{\frac{N-2}{2}}}
H(\xi_\varepsilon,\xi_\varepsilon)\frac{U_{\lambda_\varepsilon, \xi_\varepsilon}^{2_\mu^*}(y)U_{\lambda_\varepsilon, \xi_\varepsilon}^{2_\mu^*-1}(x)}{|x-y|^{\mu}}dxdy
\nonumber\\
=&-\frac{N(N-2)^2\omega_{N}B_N}{\mathcal{A}_{H,L}}\frac{H(0,0)}{\lambda_\varepsilon^{N-2}}
+O\left(\frac{1}{\lambda_\varepsilon^{N}}
+\varepsilon^N\lambda_\varepsilon^2\right)
=-\frac{N(N-2)^2\omega_{N}B_N}{\mathcal{A}_{H,L}}\frac{H(0,0)}{\lambda_\varepsilon^{N-2}}
+O\left(\varepsilon^{\frac{N}{2}}+\varepsilon^{N-1}
\right).
\end{align}
Moreover,
\begin{align}\label{B2}
&-\int_{\Omega_\varepsilon}\int_{\Omega_\varepsilon}\frac{(N-2)\omega_{N}}{\lambda_\varepsilon^{\frac{N-2}{2}}}\left\langle \nabla H(\xi_\varepsilon,\xi_\varepsilon), x-\xi_\varepsilon\right\rangle\frac{U_{\lambda_\varepsilon, \xi_\varepsilon}^{2_\mu^*}(y)U_{\lambda_\varepsilon, \xi_\varepsilon}^{2_\mu^*-1}(x)}{|x-y|^{\mu}}dxdy
\nonumber\\
=&O\left(\frac{1}{\lambda_\varepsilon^{\frac{N-2}{2}}}\int_{\Omega_\varepsilon}\int_{\Omega_\varepsilon}|\nabla H(\xi_\varepsilon,\xi_\varepsilon)||x-\xi_\varepsilon|\frac{U_{\lambda_\varepsilon, \xi_\varepsilon}^{2_\mu^*}(y)U_{\lambda_\varepsilon, \xi_\varepsilon}^{2_\mu^*-1}(x)}{|x-y|^{\mu}}dxdy
\right)
\nonumber\\
=&O\left(\frac{1}{\lambda_\varepsilon^{\frac{N-2}{2}}}\frac{N(N-2)}{\mathcal{A}_{H,L}}\int_{\Omega_\varepsilon}U_{\lambda_\varepsilon, \xi_\varepsilon}^{2^*-1}(x)|x-\xi_\varepsilon|dx
\right)
=O\left(\frac{1}{\lambda_\varepsilon^{N-1}}\right)
=O\left(\varepsilon^{\frac{N-1}{2}}
\right),
\end{align}
and
\begin{align}\label{B3}
&-\int_{\Omega_\varepsilon}\int_{\Omega_\varepsilon}\frac{(N-2)\omega_{N}}{\lambda_\varepsilon^{\frac{N-2}{2}}}
O\left(|x-\xi_\varepsilon|^2\right)
\frac{U_{\lambda_\varepsilon, \xi_\varepsilon}^{2_\mu^*}(y)U_{\lambda_\varepsilon, \xi_\varepsilon}^{2_\mu^*-1}(x)}{|x-y|^{\mu}}dxdy
\nonumber\\
=&O\left(\frac{1}{\lambda_\varepsilon^{\frac{N-2}{2}}}\frac{N(N-2)}{\mathcal{A}_{H,L}}\int_{\Omega_\varepsilon}U_{\lambda_\varepsilon, \xi_\varepsilon}^{2^*-1}(x)|x-\xi_\varepsilon|^2dx
\right)
=O\left(\frac{\ln\lambda_\varepsilon}{\lambda_\varepsilon^{N}}
\right)
=O\left(\varepsilon^{\frac{N}{2}}|\ln \varepsilon|
\right).
\end{align}
One has
\begin{align}\label{B4}
&-\frac{\varepsilon^{N-2}\lambda_\varepsilon^{\frac{N-2}{2}}}{(1+|\tau|^2)^{\frac{N-2}{2}}}
\int_{\Omega_\varepsilon}\int_{\Omega_\varepsilon}\frac{1}{|x|^{N-2}}
\frac{U_{\lambda_\varepsilon, \xi_\varepsilon}^{2_\mu^*}(y)U_{\lambda_\varepsilon, \xi_\varepsilon}^{2_\mu^*-1}(x)}{|x-y|^{\mu}}dxdy
\nonumber\\
=&-\frac{\varepsilon^{N-2}\lambda_\varepsilon^{\frac{N-2}{2}}}{(1+|\tau|^2)^{\frac{N-2}{2}}}
\Bigg\{\int_{\mathbb{R}^N}\int_{\mathbb{R}^N}
-2\int_{\mathbb{R}^N\setminus \Omega_\varepsilon}\int_{\Omega_\varepsilon}
-\int_{\mathbb{R}^N\setminus \Omega_\varepsilon}\int_{\mathbb{R}^N\setminus \Omega_\varepsilon}\Bigg\}
\frac{U_{\lambda_\varepsilon, \xi_\varepsilon}^{2_\mu^*}(y)U_{\lambda_\varepsilon, \xi_\varepsilon}^{2_\mu^*-1}(x)}{|x-y|^{\mu}}\frac{1}{|x|^{N-2}}dxdy.
\end{align}
We set $x=\lambda_\varepsilon^{-1}z$, then
\begin{align}\label{B41}
&-\frac{\varepsilon^{N-2}\lambda_\varepsilon^{\frac{N-2}{2}}}{(1+|\tau|^2)^{\frac{N-2}{2}}}
\int_{\mathbb{R}^N}\int_{\mathbb{R}^N}\frac{U_{\lambda_\varepsilon, \xi_\varepsilon}^{2_\mu^*}(y)U_{\lambda_\varepsilon, \xi_\varepsilon}^{2_\mu^*-1}(x)}{|x-y|^{\mu}}
\frac{1}{|x|^{N-2}}dxdy
\nonumber\\
=&-\frac{N(N-2)\varepsilon^{N-2}\lambda_\varepsilon^{\frac{N-2}{2}}}{\mathcal{A}_{H,L}(1+|\tau|^2)^{\frac{N-2}{2}}}
\int_{\mathbb{R}^N}\frac{U_{\lambda_\varepsilon, \xi_\varepsilon}^{2^*-1}(x)}{|x|^{N-2}}dx
=-\frac{N(N-2)\varepsilon^{N-2}\lambda_\varepsilon^{N-2}}{\mathcal{A}_{H,L}(1+|\tau|^2)^{\frac{N-2}{2}}}
\int_{\mathbb{R}^N}\frac{1}{|z|^{N-2}(1+|z-\tau|^2)^{\frac{N+2}{2}}}dz,
\end{align}
and
\begin{align}\label{B42}
&\frac{2\varepsilon^{N-2}\lambda_\varepsilon^{\frac{N-2}{2}}}{(1+|\tau|^2)^{\frac{N-2}{2}}}
\int_{\mathbb{R}^N\setminus \Omega_\varepsilon}\int_{\Omega_\varepsilon}\frac{U_{\lambda_\varepsilon, \xi_\varepsilon}^{2_\mu^*}(y)U_{\lambda_\varepsilon, \xi_\varepsilon}^{2_\mu^*-1}(x)}{|x-y|^{\mu}}\frac{1}{|x|^{N-2}}dxdy
\nonumber\\
=&O\left(\varepsilon^{N-2}\lambda_\varepsilon^{\frac{N-2}{2}}\int_{\mathbb{R}^N\setminus \Omega_\varepsilon}
\frac{U_{\lambda_\varepsilon, \xi_\varepsilon}^{2^*-1}(x)}{|x|^{N-2}}dx
\right)
=O\left(\left(\varepsilon\lambda_\varepsilon\right)^{N-2}\Bigg[\left(\varepsilon\lambda_\varepsilon\right)^2
+\frac{1}{\lambda_\varepsilon^N}\Bigg]\right)
=O\left(\varepsilon^{\frac{N-2}{2}}\left(
\varepsilon+\varepsilon^{\frac{N}{2}}
\right)
\right).
\end{align}
We also can get
\begin{align}\label{B43}
\frac{\varepsilon^{N-2}\lambda_\varepsilon^{\frac{N-2}{2}}}{(1+|\tau|^2)^{\frac{N-2}{2}}}
\int_{\mathbb{R}^N\setminus \Omega_\varepsilon}\int_{\mathbb{R}^N\setminus \Omega_\varepsilon}\frac{U_{\lambda_\varepsilon, \xi_\varepsilon}^{2_\mu^*}(y)U_{\lambda_\varepsilon, \xi_\varepsilon}^{2_\mu^*-1}(x)}{|x-y|^{\mu}}
\frac{1}{|x|^{N-2}}dxdy
=&O\left(\left(\varepsilon\lambda_\varepsilon\right)^{N-2}\Bigg[\left(\varepsilon\lambda_\varepsilon\right)^2
+\frac{1}{\lambda_\varepsilon^N}\Bigg]\right)
\nonumber\\
=&O\left(\varepsilon^{\frac{N-2}{2}}\left(
\varepsilon+\varepsilon^{\frac{N}{2}}
\right)
\right).
\end{align}
In view of (\ref{B4}), (\ref{B41})-(\ref{B43}),
\begin{align}\label{B4h}
&-\frac{\varepsilon^{N-2}\lambda_\varepsilon^{\frac{N-2}{2}}}{(1+|\tau|^2)^{\frac{N-2}{2}}}
\int_{\Omega_\varepsilon}\int_{\Omega_\varepsilon}\frac{1}{|x|^{N-2}}
\frac{U_{\lambda_\varepsilon, \xi_\varepsilon}^{2_\mu^*}(y)U_{\lambda_\varepsilon, \xi_\varepsilon}^{2_\mu^*-1}(x)}{|x-y|^{\mu}}dxdy
\nonumber\\
=&-\frac{N(N-2)}{\mathcal{A}_{H,L}(1+|\tau|^2)^{\frac{N-2}{2}}}
\left(\int_{\mathbb{R}^N}\frac{1}{|z|^{N-2}(1+|z-\tau|^2)^{\frac{N+2}{2}}}dz\right)
\varepsilon^{N-2}\lambda_\varepsilon^{N-2}\left(1+o(1)\right).
\end{align}
It follows from Lemma \ref{PU} that
\begin{align}\label{B5}
\int_{\Omega_\varepsilon}\int_{\Omega_\varepsilon}\frac{U_{\lambda_\varepsilon, \xi_\varepsilon}^{2_\mu^*}(y)
U_{\lambda_\varepsilon, \xi_\varepsilon}^{2_\mu^*-1}(x)R(x)}{|x-y|^{\mu}}dxdy
=O\left(\varepsilon^{N-2}
+\varepsilon^{N-1}\lambda_\varepsilon^{N-1}
+\frac{1}{\lambda_\varepsilon^N}\right)
=O\left(\varepsilon^{N-2}
+\varepsilon^{\frac{N-1}{2}}\right).
\end{align}
Then combining (\ref{B}) with (\ref{B1}), (\ref{B2}), (\ref{B3}), (\ref{B4h}), (\ref{B5}), we have
\begin{align}\label{bh}
B=&-\frac{N(N-2)^2\omega_{N}B_N}{\mathcal{A}_{H,L}}\frac{H(0,0)}{\lambda_\varepsilon^{N-2}}\left(1+o(1)\right)
\nonumber\\
&-\frac{N(N-2)}{\mathcal{A}_{H,L}(1+|\tau|^2)^{\frac{N-2}{2}}}
\left(\int_{\mathbb{R}^N}\frac{1}{|z|^{N-2}(1+|z-\tau|^2)^{\frac{N+2}{2}}}dz\right)
\varepsilon^{N-2}\lambda_\varepsilon^{N-2}\left(1+o(1)\right).
\end{align}
Finally, we compute $C$ in (\ref{IPU1}). It holds:
\begin{align}\label{C}
C=&\Bigg\{\int_{\mathbb{R}^N}\int_{\mathbb{R}^N}
-2\int_{\mathbb{R}^N\setminus \Omega_\varepsilon}\int_{\Omega_\varepsilon}
-\int_{\mathbb{R}^N\setminus \Omega_\varepsilon}\int_{\mathbb{R}^N\setminus \Omega_\varepsilon}\Bigg\}
\frac{PU_{\lambda_\varepsilon, \xi_\varepsilon}^{2_\mu^*}(y)PU_{\lambda_\varepsilon, \xi_\varepsilon}^{2_\mu^*}(x)}{|x-y|^{\mu}}dxdy
\nonumber\\
:=&C_0-2C_1-C_2.
\end{align}
Since
\begin{align*}
PU_{\lambda_\varepsilon, \xi_\varepsilon}^{2_\mu^*}=U_{\lambda_\varepsilon, \xi_\varepsilon}^{2_\mu^*}-2_\mu^*U_{\lambda_\varepsilon, \xi_\varepsilon}^{2_\mu^*-1}\psi_{\lambda_\varepsilon, \xi_\varepsilon}
+O\left(U_{\lambda_\varepsilon, \xi_\varepsilon}^{2_\mu^*-2}\psi^2_{\lambda_\varepsilon, \xi_\varepsilon}\right),
\end{align*}
where $\psi_{\lambda_\varepsilon, \xi_\varepsilon}=U_{\lambda_\varepsilon, \xi_\varepsilon}-PU_{\lambda_\varepsilon, \xi_\varepsilon}$. We have
\begin{align}\label{C0}
C_0=&\int_{\mathbb{R}^N}\int_{\mathbb{R}^N}\frac{U_{\lambda_\varepsilon, \xi_\varepsilon}^{2_\mu^*}(y)U_{\lambda_\varepsilon, \xi_\varepsilon}^{2_\mu^*}(x)}{|x-y|^{\mu}}dxdy
-2\cdot 2_\mu^*\int_{\mathbb{R}^N}\int_{\mathbb{R}^N}\frac{U_{\lambda_\varepsilon, \xi_\varepsilon}^{2_\mu^*}(y)U_{\lambda_\varepsilon, \xi_\varepsilon}^{2_\mu^*-1}(x)\psi_{\lambda_\varepsilon, \xi_\varepsilon}(x)}{|x-y|^{\mu}}dxdy
\nonumber\\
&+O\left(C_{01}+C_{02}+C_{03}+C_{04}\right),
\end{align}
where
\begin{align}\label{C01}
C_{01}=\int_{\mathbb{R}^N}\int_{\mathbb{R}^N}\frac{U_{\lambda_\varepsilon, \xi_\varepsilon}^{2_\mu^*}(y)U_{\lambda_\varepsilon, \xi_\varepsilon}^{2_\mu^*-2}(x)\psi^2_{\lambda_\varepsilon, \xi_\varepsilon}(x)}{|x-y|^{\mu}}dxdy,
\end{align}
\begin{align}\label{C02}
C_{02}=\int_{\mathbb{R}^N}\int_{\mathbb{R}^N}\frac{U_{\lambda_\varepsilon, \xi_\varepsilon}^{2_\mu^*-1}(y)\psi_{\lambda_\varepsilon, \xi_\varepsilon}(y)U_{\lambda_\varepsilon, \xi_\varepsilon}^{2_\mu^*-1}(x)\psi_{\lambda_\varepsilon, \xi_\varepsilon}(x)}{|x-y|^{\mu}}dxdy,
\end{align}
\begin{align}\label{C03}
C_{03}=\int_{\mathbb{R}^N}\int_{\mathbb{R}^N}\frac{U_{\lambda_\varepsilon, \xi_\varepsilon}^{2_\mu^*-1}(y)\psi_{\lambda_\varepsilon, \xi_\varepsilon}(y)U_{\lambda_\varepsilon, \xi_\varepsilon}^{2_\mu^*-2}(x)\psi^2_{\lambda_\varepsilon, \xi_\varepsilon}(x)}{|x-y|^{\mu}}dxdy,
\end{align}
\begin{align}\label{C04}
C_{04}=\int_{\mathbb{R}^N}\int_{\mathbb{R}^N}\frac{U_{\lambda_\varepsilon, \xi_\varepsilon}^{2_\mu^*-2}(y)\psi^2_{\lambda_\varepsilon, \xi_\varepsilon}(y)U_{\lambda_\varepsilon, \xi_\varepsilon}^{2_\mu^*-2}(x)\psi^2_{\lambda_\varepsilon, \xi_\varepsilon}(x)}{|x-y|^{\mu}}dxdy.
\end{align}
We have
\begin{align}\label{C010}
C_{01}
=\frac{N(N-2)}{\mathcal{A}_{H,L}}\int_{\mathbb{R}^N}U_{\lambda_\varepsilon, \xi_\varepsilon}^{2^*-2}(x)\psi^2_{\lambda_\varepsilon, \xi_\varepsilon}(x)dx
=\frac{N(N-2)}{\mathcal{A}_{H,L}}\left(\int_{\Omega_\varepsilon}+\int_{\mathbb{R}^N\setminus\Omega_\varepsilon}\right)
U_{\lambda_\varepsilon, \xi_\varepsilon}^{2^*-2}(x)\psi^2_{\lambda_\varepsilon, \xi_\varepsilon}(x)dx.
\end{align}
By \cite[Lemma 3.2]{gemp}, we deduce that
\begin{align*}
\int_{\Omega_\varepsilon}U_{\lambda_\varepsilon, \xi_\varepsilon}^{2^*-2}(x)\psi^2_{\lambda_\varepsilon, \xi_\varepsilon}(x)dx
=O\left(\frac{1}{\lambda_\varepsilon^N}
+(\varepsilon\lambda_\varepsilon)^N\right)
=O\left(\varepsilon^{\frac{N}{2}}\right).
\end{align*}
It follows from $0<\psi_{\lambda_\varepsilon, \xi_\varepsilon}\leq U_{\lambda_\varepsilon, \xi_\varepsilon}$ that
\begin{align*}
\int_{\mathbb{R}^N\setminus\Omega_\varepsilon}U_{\lambda_\varepsilon, \xi_\varepsilon}^{2^*-2}(x)\psi^2_{\lambda_\varepsilon, \xi_\varepsilon}(x)dx
=O\left(\int_{\mathbb{R}^N\setminus\Omega_\varepsilon}U_{\lambda_\varepsilon, \xi_\varepsilon}^{2^*}(x)dx\right)
=O\left(\frac{1}{\lambda_\varepsilon^N}+(\varepsilon\lambda_\varepsilon)^N\right)
=O\left(\varepsilon^{\frac{N}{2}}\right).
\end{align*}
Above all,
\begin{align}\label{C01h}
C_{01}=O\left(\frac{1}{\lambda_\varepsilon^N}
+(\varepsilon\lambda_\varepsilon)^N\right)
=O\left(\varepsilon^{\frac{N}{2}}\right).
\end{align}
Using Hardy-Littlewood-Sobolev inequality and Lemma \ref{PU},
\begin{align}\label{C020}
C_{02}=&\int_{\mathbb{R}^N}\int_{\mathbb{R}^N}\frac{U_{\lambda_\varepsilon, \xi_\varepsilon}^{2_\mu^*-1}(y)\psi_{\lambda_\varepsilon, \xi_\varepsilon}(y)U_{\lambda_\varepsilon, \xi_\varepsilon}^{2_\mu^*-1}(x)\psi_{\lambda_\varepsilon, \xi_\varepsilon}(x)}{|x-y|^{\mu}}dxdy
\nonumber\\
\leq &c\left(\int_{\mathbb{R}^N}\left(U_{\lambda_\varepsilon, \xi_\varepsilon}^{2_\mu^*-1}(x)\psi_{\lambda_\varepsilon, \xi_\varepsilon}(x)\right)^{\frac{2N}{2N-\mu}}dx\right)^{\frac{2(2N-\mu)}{2N}}
\nonumber\\
\leq &c\left(\left(\int_{\Omega_\varepsilon}+\int_{\mathbb{R}^N\setminus\Omega_\varepsilon}\right)\left(U_{\lambda_\varepsilon, \xi_\varepsilon}^{2_\mu^*-1}(x)\psi_{\lambda_\varepsilon, \xi_\varepsilon}(x)\right)^{\frac{2N}{2N-\mu}}dx\right)^{\frac{2(2N-\mu)}{2N}}
\nonumber\\
\leq &c\left(\int_{\Omega_\varepsilon}\left(U_{\lambda_\varepsilon, \xi_\varepsilon}^{2_\mu^*-1}(x)\psi_{\lambda_\varepsilon, \xi_\varepsilon}(x)\right)^{\frac{2N}{2N-\mu}}dx
+\int_{\mathbb{R}^N\setminus\Omega_\varepsilon}
\left(U_{\lambda_\varepsilon, \xi_\varepsilon}^{2_\mu^*}(x)\right)^{\frac{2N}{2N-\mu}}dx\right)^{\frac{2(2N-\mu)}{2N}}
\nonumber\\
=&O\left(\frac{1}{\lambda_\varepsilon^{2N-4}}+\varepsilon^{2N-4}\lambda_\varepsilon^{2N-4}
+\varepsilon^{2N-\mu}\lambda_\varepsilon^{2N-\mu}+\frac{1}{\lambda_\varepsilon^{2N-\mu}}\right)
=o\left(\frac{1}{\lambda_\varepsilon^{N-2}}
+(\varepsilon\lambda_\varepsilon)^{N-2}\right)
=o\left(\varepsilon^{\frac{N-2}{2}}\right).
\end{align}
And similar to the estimate of $C_{02}$, we can also have
\begin{align}\label{C034}
C_{03}=C_{04}
=o\left(\frac{1}{\lambda_\varepsilon^{N-2}}
+(\varepsilon\lambda_\varepsilon)^{N-2}\right)
=o\left(\varepsilon^{\frac{N-2}{2}}\right).
\end{align}
For the second term in (\ref{C0}), from the estimate of $B$, we get
\begin{align}\label{C0second}
&-2\cdot 2_\mu^*\int_{\mathbb{R}^N}\int_{\mathbb{R}^N}\frac{U_{\lambda_\varepsilon, \xi_\varepsilon}^{2_\mu^*}(y)
U_{\lambda_\varepsilon, \xi_\varepsilon}^{2_\mu^*-1}(x)\psi_{\lambda_\varepsilon, \xi_\varepsilon}(x)}{|x-y|^{\mu}}dxdy
\nonumber\\
=&2\cdot 2_\mu^*\Bigg\{-\frac{N(N-2)^2\omega_{N}B_N}{\mathcal{A}_{H,L}}\frac{H(0,0)}{\lambda_\varepsilon^{N-2}}\left(1+o(1)\right)
\nonumber\\
&-\frac{N(N-2)}{\mathcal{A}_{H,L}(1+|\tau|^2)^{\frac{N-2}{2}}}
\left(\int_{\mathbb{R}^N}\frac{1}{|z|^{N-2}(1+|z-\tau|^2)^{\frac{N+2}{2}}}dz\right)
\varepsilon^{N-2}\lambda_\varepsilon^{N-2}\left(1+o(1)\right)
\Bigg\}.
\end{align}
Then by the estimate of $A$ and (\ref{C0}), (\ref{C01h}), (\ref{C020}), (\ref{C034}), (\ref{C0second}), we obtain
\begin{align}\label{C0h}
C_0=&\frac{N(N-2)}{\mathcal{A}_{H,L}}A_N-2\cdot 2_\mu^*\Bigg\{\frac{N(N-2)^2\omega_{N}B_N}{\mathcal{A}_{H,L}}\frac{H(0,0)}{\lambda_\varepsilon^{N-2}}\left(1+o(1)\right)
\nonumber\\
&+\frac{N(N-2)}{\mathcal{A}_{H,L}(1+|\tau|^2)^{\frac{N-2}{2}}}
\left(\int_{\mathbb{R}^N}\frac{1}{|z|^{N-2}(1+|z-\tau|^2)^{\frac{N+2}{2}}}dz\right)
\varepsilon^{N-2}\lambda_\varepsilon^{N-2}\left(1+o(1)\right)
\Bigg\}.
\end{align}
Similar to $B_1$ and $B_2$, we have
\begin{align}\label{C12}
C_1=C_2
=O\left(\frac{1}{\lambda_\varepsilon^N}
+(\varepsilon\lambda_\varepsilon)^N\right)
=O\left(\varepsilon^{\frac{N}{2}}\right).
\end{align}
Combing (\ref{C}) with (\ref{C0h}), (\ref{C12}), we can prove
\begin{align}\label{Ch}
C=&\frac{N(N-2)}{\mathcal{A}_{H,L}}A_N-2\cdot 2_\mu^*\Bigg\{\frac{N(N-2)^2\omega_{N}B_N}{\mathcal{A}_{H,L}}\frac{H(0,0)}{\lambda_\varepsilon^{N-2}}\left(1+o(1)\right)
\nonumber\\
&+\frac{N(N-2)}{\mathcal{A}_{H,L}(1+|\tau|^2)^{\frac{N-2}{2}}}
\left(\int_{\mathbb{R}^N}\frac{1}{|z|^{N-2}(1+|z-\tau|^2)^{\frac{N+2}{2}}}dz\right)
\varepsilon^{N-2}\lambda_\varepsilon^{N-2}\left(1+o(1)\right)
\Bigg\}.
\end{align}
Consequently, it follows from  (\ref{IPU1}), (\ref{ah}), (\ref{bh}) and (\ref{Ch}) that
\begin{align*}
I_\varepsilon(PU_{\lambda, \xi})
=&\frac{1}{2}A+\frac{1}{2}B-\frac{1}{2\cdot2_\mu^*}C
\nonumber\\
=&\frac{1}{2}\Bigg\{\frac{N(N-2)}{\mathcal{A}_{H,L}}A_N
+O\left(\varepsilon^{\frac{N}{2}}\right)
\Bigg\}
+\frac{1}{2}\Bigg\{-\frac{N(N-2)^2\omega_{N}B_N}{\mathcal{A}_{H,L}}\frac{H(0,0)}{\lambda_\varepsilon^{N-2}}\left(1+o(1)\right)
\nonumber\\
&-\frac{N(N-2)}{\mathcal{A}_{H,L}(1+|\tau|^2)^{\frac{N-2}{2}}}
\left(\int_{\mathbb{R}^N}\frac{1}{|z|^{N-2}(1+|z-\tau|^2)^{\frac{N+2}{2}}}dz\right)
\varepsilon^{N-2}\lambda_\varepsilon^{N-2}\left(1+o(1)\right)
\Bigg\}
\nonumber\\
&-\frac{1}{2\cdot2_\mu^*}\Bigg\{\frac{N(N-2)}{\mathcal{A}_{H,L}}A_N-2\cdot 2_\mu^*\Bigg\{\frac{N(N-2)^2\omega_{N}B_N}{\mathcal{A}_{H,L}}\frac{H(0,0)}{\lambda_\varepsilon^{N-2}}\left(1+o(1)\right)
\nonumber\\
&+\frac{N(N-2)}{\mathcal{A}_{H,L}(1+|\tau|^2)^{\frac{N-2}{2}}}
\left(\int_{\mathbb{R}^N}\frac{1}{|z|^{N-2}(1+|z-\tau|^2)^{\frac{N+2}{2}}}dz\right)
\varepsilon^{N-2}\lambda_\varepsilon^{N-2}\left(1+o(1)\right)
\Bigg\}
\Bigg\}
\nonumber\\
=&\frac{1}{2}\left(1-\frac{1}{2_\mu^*}\right)\frac{N(N-2)}{\mathcal{A}_{H,L}}A_N
+\frac{N(N-2)^2\omega_{N}B_N}{2\mathcal{A}_{H,L}}\frac{H(0,0)}{\lambda_\varepsilon^{N-2}}\left(1+o(1)\right)
\nonumber\\
&+\frac{N(N-2)}{2\mathcal{A}_{H,L}(1+|\tau|^2)^{\frac{N-2}{2}}}
\left(\int_{\mathbb{R}^N}\frac{1}{|z|^{N-2}(1+|z-\tau|^2)^{\frac{N+2}{2}}}dz\right)
\varepsilon^{N-2}\lambda_\varepsilon^{N-2}\left(1+o(1)\right).
\end{align*}
Next, we give the $C^1$ estimate of $I_\varepsilon(PU_{\lambda_\varepsilon, \xi_\varepsilon})$. Let $\partial_s$ denote $\partial_{\lambda}$ and $\partial_{\tau^i}$ for $i=1,\cdots,N$. First of all, let $\partial_s=\partial_{\lambda}$ and $PU_{\lambda_{\varepsilon}, \xi_{\varepsilon}}:=PU$, we have
\begin{align}\label{partialI}
&\partial_sI_\varepsilon(PU)
\nonumber\\
=&\int_{\Omega_\varepsilon}\nabla PU\nabla\left(\partial_sPU\right)
-\int_{\Omega_\varepsilon}\int_{\Omega_\varepsilon}\frac{PU^{2_\mu^*}(y)PU^{2_\mu^*-1}(x)\partial_sPU(x)}{|x-y|^{\mu}}dxdy
\nonumber\\
=&\int_{\Omega_\varepsilon}\int_{\Omega_\varepsilon}\frac{U^{2_\mu^*}(y)U^{2_\mu^*-1}(x)\partial_sPU(x)}{|x-y|^{\mu}}dxdy
-\int_{\Omega_\varepsilon}\int_{\Omega_\varepsilon}\frac{PU^{2_\mu^*}(y)PU^{2_\mu^*-1}(x)\partial_sPU(x)}{|x-y|^{\mu}}dxdy
\nonumber\\
=&\int_{\Omega_\varepsilon}\int_{\Omega_\varepsilon}\frac{U^{2_\mu^*}(y)U^{2_\mu^*-1}(x)\partial_sU(x)}{|x-y|^{\mu}}dxdy
+\int_{\Omega_\varepsilon}\int_{\Omega_\varepsilon}\frac{U^{2_\mu^*}(y)U^{2_\mu^*-1}(x)
\left(\partial_sPU(x)-\partial_sU(x)\right)
}{|x-y|^{\mu}}dxdy
\nonumber\\
&
-\int_{\Omega_\varepsilon}\int_{\Omega_\varepsilon}\frac{PU^{2_\mu^*}(y)PU^{2_\mu^*-1}(x)\partial_sPU(x)}{|x-y|^{\mu}}dxdy
\nonumber\\
:=&A'+B'-C'.
\end{align}
We next estimate $A'$, $B'$, $C'$ respectively. It holds
\begin{align*}
A'=&\Bigg\{\int_{\mathbb{R}^N}\int_{\mathbb{R}^N}
-2\int_{\mathbb{R}^N\setminus \Omega_\varepsilon}\int_{\Omega_\varepsilon}
-\int_{\mathbb{R}^N\setminus \Omega_\varepsilon}\int_{\mathbb{R}^N\setminus \Omega_\varepsilon}
\Bigg\}
\frac{U^{2_\mu^*}(y)U^{2_\mu^*-1}(x)\partial_sU(x)}{|x-y|^{\mu}}dxdy
\nonumber\\
:=&A'_0-2A'_1-A'_2,
\end{align*}
where
\begin{align}\label{A'0}
A'_0=\frac{N(N-2)}{\mathcal{A}_{H,L}}\int_{\mathbb{R}^N}
U^{2^*-1}(x)\partial_sU(x)dx
=\frac{N(N-2)}{2^*\mathcal{A}_{H,L}}\int_{\mathbb{R}^N}\frac{\partial}{\partial s}\left(U^{2^*}(x)\right)dx
=0,
\end{align}
and since $\frac{\partial U}{\partial \lambda_{\varepsilon}}=O\left(\frac{U}{\lambda_{\varepsilon}}\right)$ and $\frac{\partial \lambda_{\varepsilon}}{\partial \lambda}=\varepsilon^{-\frac{1}{2}}:=\varepsilon^{-\theta}$, we have
\begin{align}\label{A'1}
A'_1=O\left(\int_{\mathbb{R}^N\setminus \Omega_\varepsilon}\int_{\Omega_\varepsilon}
\frac{U^{2_\mu^*}(y)U^{2_\mu^*}(x)}{|x-y|^{\mu}}dxdy
\right)
=O\left(A_1\right)
=O\left(\frac{1}{\lambda_{\varepsilon}^{N}}
+\varepsilon^N\lambda_{\varepsilon}^{N}\right)
=O\left(\varepsilon^{\frac{N}{2}}\right).
\end{align}
Similarly, we can also get
\begin{align}\label{A'2}
A'_2=O\left(\frac{1}{\lambda_{\varepsilon}^{N}}
+\varepsilon^N\lambda_{\varepsilon}^{N}\right)
=O\left(\varepsilon^{\frac{N}{2}}\right).
\end{align}
By (\ref{A'0}), (\ref{A'1}), (\ref{A'2}), we conclude that
\begin{align}\label{A'}
A'=O\left(\frac{1}{\lambda_{\varepsilon}^{N}}
+\varepsilon^N\lambda_{\varepsilon}^{N}\right)
=O\left(\varepsilon^{\frac{N}{2}}\right).
\end{align}
We next estimate $B'$. In terms of Lemma \ref{PU}, we can get
\begin{align*}
B'=&\int_{\Omega_\varepsilon}\int_{\Omega_\varepsilon}\frac{U^{2_\mu^*}(y)U^{2_\mu^*-1}(x)
\partial_s\left(PU(x)-U(x)\right)
}{|x-y|^{\mu}}dxdy
\nonumber\\
=&\int_{\Omega_\varepsilon}\int_{\Omega_\varepsilon}\frac{U^{2_\mu^*}(y)U^{2_\mu^*-1}(x)
\left(\frac{(N-2)^2\omega_{N}}{2\lambda_{\varepsilon}^{\frac{N}{2}}}H(x,\xi_{\varepsilon})\varepsilon^{-\theta}
-\frac{N-2}{2}\frac{\lambda_{\varepsilon}^{\frac{N-4}{2}}\varepsilon^{-\theta}}{(1+|\tau|^2)^{\frac{N-2}{2}}}\frac{\varepsilon^{N-2}}{|x|^{N-2}}
+\partial_sR(x)
\right)
}{|x-y|^{\mu}}dxdy.
\end{align*}
By Taylor expansion, we have
\begin{align}\label{B'0}
&\frac{(N-2)^2\omega_{N}}{2\lambda_{\varepsilon}^{\frac{N}{2}}}\int_{\Omega_\varepsilon}\int_{\Omega_\varepsilon}
\frac{U^{2_\mu^*}(y)U^{2_\mu^*-1}(x)
H(x,\xi_{\varepsilon})\varepsilon^{-\theta}}{|x-y|^{\mu}}dxdy
\nonumber\\
=&\frac{(N-2)^2\omega_{N}}{2\lambda_{\varepsilon}^{\frac{N}{2}}}\int_{\Omega_\varepsilon}\int_{\Omega_\varepsilon}
\frac{U^{2_\mu^*}(y)U^{2_\mu^*-1}(x)
\left(H(\xi_{\varepsilon},\xi_{\varepsilon})+\left\langle \nabla H(\xi_{\varepsilon},\xi_{\varepsilon}), x-\xi_{\varepsilon}\right\rangle+O\left(|x-\xi_{\varepsilon}|^2\right)
\right)
\varepsilon^{-\theta}}{|x-y|^{\mu}}dxdy
\nonumber\\
=&\frac{N(N-2)^3\omega_{N}B_{N}}{2\mathcal{A}_{H,L}\lambda_{\varepsilon}^{N-1}}H(0,0)\varepsilon^{-\theta}
+O\left(\varepsilon^N\lambda_{\varepsilon}^2\right)
+O\left(\frac{1}{\lambda_{\varepsilon}^{N-1}}\right)
+O\left(\frac{\ln \lambda_{\varepsilon}}{\lambda_{\varepsilon}^{N}}\right)
\nonumber\\
=&\frac{N(N-2)^3\omega_{N}B_{N}}{2\mathcal{A}_{H,L}\lambda_{\varepsilon}^{N-1}}H(0,0)\varepsilon^{-\theta}
+O\left(\varepsilon^{N-1}\right)
+O\left(\varepsilon^{\frac{N-1}{2}}\right)
+O\left(\varepsilon^{\frac{N}{2}}|\ln \varepsilon|\right).
\end{align}
Moreover,
\begin{align}\label{B'1}
&-\frac{(N-2)}{2}\int_{\Omega_\varepsilon}\int_{\Omega_\varepsilon}\frac{U^{2_\mu^*}(y)U^{2_\mu^*-1}(x)
\frac{\lambda_{\varepsilon}^{\frac{N-4}{2}}\varepsilon^{-\theta}}{(1+|\tau|^2)^{\frac{N-2}{2}}}\frac{\varepsilon^{N-2}}{|x|^{N-2}}
}{|x-y|^{\mu}}dxdy
\nonumber\\
=&-\frac{(N-2)}{2}\Bigg\{\int_{\mathbb{R}^N}\int_{\mathbb{R}^N}
-2\int_{\mathbb{R}^N\setminus \Omega_\varepsilon}\int_{\Omega_\varepsilon}
-\int_{\mathbb{R}^N\setminus \Omega_\varepsilon}\int_{\mathbb{R}^N\setminus \Omega_\varepsilon}
\Bigg\}
\frac{U^{2_\mu^*}(y)U^{2_\mu^*-1}(x)
\frac{\lambda_{\varepsilon}^{\frac{N-4}{2}}\varepsilon^{-\theta}}{(1+|\tau|^2)^{\frac{N-2}{2}}}\frac{\varepsilon^{N-2}}{|x|^{N-2}}
}{|x-y|^{\mu}}dxdy
\nonumber\\
=&-\frac{N(N-2)^2}{2\mathcal{A}_{H,L}}\frac{M(\tau)}{(1+|\tau|^2)^{\frac{N-2}{2}}}
\varepsilon^{N-2}\lambda_{\varepsilon}^{N-3}\varepsilon^{-\theta}
+O\left(\varepsilon^{N-2}\lambda_{\varepsilon}^{N-2}\left((\varepsilon\lambda_{\varepsilon})^2
+\frac{1}{\lambda_{\varepsilon}^N}\right)\right)
\nonumber\\
=&-\frac{N(N-2)^2}{2\mathcal{A}_{H,L}}\frac{M(\tau)}{(1+|\tau|^2)^{\frac{N-2}{2}}}
\varepsilon^{N-2}\lambda_{\varepsilon}^{N-3}\varepsilon^{-\theta}
+O\left(\varepsilon^{\frac{N-2}{2}}\left(
\varepsilon+\varepsilon^{\frac{N}{2}}\right)\right).
\end{align}
By Lemma \ref{PU}, we have
\begin{align}\label{B'2}
\int_{\Omega_\varepsilon}\int_{\Omega_\varepsilon}\frac{U^{2_\mu^*}(y)U^{2_\mu^*-1}(x)\partial_sR(x)
}{|x-y|^{\mu}}dxdy
=o\left(\lambda_{\varepsilon}^{N-2}\varepsilon^{N-2}
+\frac{1}{\lambda_{\varepsilon}^{N-2}}\right)
=o\left(\varepsilon^{\frac{N-2}{2}}\right).
\end{align}
Then from (\ref{B'0}), (\ref{B'1}) and (\ref{B'2}), we prove that
\begin{align}\label{B'}
B'=\frac{N(N-2)^3\omega_{N}B_{N}}{2\mathcal{A}_{H,L}\lambda_{\varepsilon}^{N-1}}H(0,0)\varepsilon^{-\theta}(1+o(1))
-\frac{N(N-2)^2}{2\mathcal{A}_{H,L}}\frac{M(\tau)}{(1+|\tau|^2)^{\frac{N-2}{2}}}
\varepsilon^{N-2}\lambda_{\varepsilon}^{N-3}\varepsilon^{-\theta}(1+o(1)).
\end{align}
As for $C'$, we have
\begin{align*}
C'=&\int_{\Omega_\varepsilon}\int_{\Omega_\varepsilon}\frac{PU^{2_\mu^*}(y)PU^{2_\mu^*-1}(x)\partial_sPU(x)}
{|x-y|^{\mu}}dxdy
\nonumber\\
=&\Bigg\{\int_{\mathbb{R}^N}\int_{\mathbb{R}^N}
-2\int_{\mathbb{R}^N\setminus \Omega_\varepsilon}\int_{\Omega_\varepsilon}
-\int_{\mathbb{R}^N\setminus \Omega_\varepsilon}\int_{\mathbb{R}^N\setminus \Omega_\varepsilon}
\Bigg\}
\frac{PU^{2_\mu^*}(y)PU^{2_\mu^*-1}(x)\partial_sPU(x)}{|x-y|^{\mu}}dxdy
\nonumber\\
=&C'_0-2C'_1-C'_2.
\end{align*}
We first estimate $C'_0$. In fact,
\begin{align*}
C'_0=&\int_{\mathbb{R}^N}\int_{\mathbb{R}^N}\frac{PU^{2_\mu^*}(y)\partial_s\left(\frac{1}{2_\mu^*}PU^{2_\mu^*}(x)\right)}{|x-y|^{\mu}}dxdy
\nonumber\\
=&\frac{1}{2_\mu^*}\int_{\mathbb{R}^N}\int_{\mathbb{R}^N}\frac{U^{2_\mu^*}(y)\partial_sU^{2_\mu^*}(x)
-2_\mu^*U^{2_\mu^*}(y)\partial_s\left(U^{2_\mu^*-1}(x)\psi(x)\right)
-2_\mu^*U^{2_\mu^*-1}(y)\psi(y)\partial_sU^{2_\mu^*}(x)
}{|x-y|^{\mu}}dxdy
\nonumber\\
&+O\Bigg(\int_{\mathbb{R}^N}\int_{\mathbb{R}^N}\frac{U^{2_\mu^*}(y)\partial_s\left(U^{2_\mu^*-2}(x)\psi^2(x)\right)}{|x-y|^{\mu}}dxdy
+\int_{\mathbb{R}^N}\int_{\mathbb{R}^N}\frac{U^{2_\mu^*-1}(y)\psi(y)\partial_s\left(U^{2_\mu^*-1}(x)\psi(x)\right)}{|x-y|^{\mu}}dxdy
\nonumber\\
&+\int_{\mathbb{R}^N}\int_{\mathbb{R}^N}\frac{U^{2_\mu^*-1}(y)\psi(y)\partial_s\left(U^{2_\mu^*-2}(x)\psi^2(x)\right)}{|x-y|^{\mu}}dxdy
+\int_{\mathbb{R}^N}\int_{\mathbb{R}^N}\frac{U^{2_\mu^*-2}(y)\psi^2(y)\partial_s\left(U^{2_\mu^*-2}(x)\psi^2(x)\right)}{|x-y|^{\mu}}dxdy
\Bigg),
\end{align*}
where
\begin{align*}
\frac{1}{2_\mu^*}\int_{\mathbb{R}^N}\int_{\mathbb{R}^N}\frac{U^{2_\mu^*}(y)\partial_sU^{2_\mu^*}(x)
}{|x-y|^{\mu}}dxdy
=\frac{N(N-2)}{\mathcal{A}_{H,L}}\int_{\mathbb{R}^N}U^{2^*-1}(x)\partial_sU(x)dx
=0,
\end{align*}
\begin{align*}
&\frac{1}{2_\mu^*}\int_{\mathbb{R}^N}\int_{\mathbb{R}^N}\frac{
-2_\mu^*U^{2_\mu^*}(y)\partial_s\left(U^{2_\mu^*-1}(x)\psi(x)\right)
-2_\mu^*U^{2_\mu^*-1}(y)\psi(y)\partial_sU^{2_\mu^*}(x)
}{|x-y|^{\mu}}dxdy
\nonumber\\
=&-\partial_s\Bigg[\int_{\mathbb{R}^N}\int_{\mathbb{R}^N}\frac{U^{2_\mu^*}(y)U^{2_\mu^*-1}(x)\psi(x)}{|x-y|^{\mu}}dxdy
\Bigg]
\nonumber\\
=&-\frac{N(N-2)}{\mathcal{A}_{H,L}}\partial_s\Bigg[\int_{\mathbb{R}^N}U^{2_\mu^*-1}(x)\psi(x)dx
\Bigg]
\nonumber\\
=&\partial_s\Bigg[-\frac{N(N-2)^2\omega_{N}B_N}{\mathcal{A}_{H,L}}\frac{H(0,0)}{\lambda_{\varepsilon}^{N-2}}\left(1+o(1)\right)
-\frac{N(N-2)M(\tau)}{\mathcal{A}_{H,L}(1+|\tau|^2)^{\frac{N-2}{2}}}
\varepsilon^{N-2}\lambda_{\varepsilon}^{N-2}\left(1+o(1)\right)
\Bigg]
\nonumber\\
=&\frac{N(N-2)^2}{\mathcal{A}_{H,L}}\left((N-2)\omega_{N}B_N\frac{H(0,0)\varepsilon^{-\theta}}{\lambda_{\varepsilon}^{N-1}}
-\frac{M(\tau)}{(1+|\tau|^2)^{\frac{N-2}{2}}}
\varepsilon^{N-2}\lambda_{\varepsilon}^{N-3}\varepsilon^{-\theta}\left(1+o(1)\right)
\right).
\end{align*}
Similar to the estimate of $C_{01}, C_{02}$, we find
\begin{align*}
&O\Bigg(\int_{\mathbb{R}^N}\int_{\mathbb{R}^N}\frac{U^{2_\mu^*}(y)\partial_s\left(U^{2_\mu^*-2}(x)\psi^2(x)\right)}{|x-y|^{\mu}}dxdy
+\int_{\mathbb{R}^N}\int_{\mathbb{R}^N}\frac{U^{2_\mu^*-1}(y)\psi(y)\partial_s\left(U^{2_\mu^*-1}(x)\psi(x)\right)}{|x-y|^{\mu}}dxdy
\nonumber\\
&+\int_{\mathbb{R}^N}\int_{\mathbb{R}^N}\frac{U^{2_\mu^*-1}(y)\psi(y)\partial_s\left(U^{2_\mu^*-2}(x)\psi^2(x)\right)}{|x-y|^{\mu}}dxdy
\nonumber\\
&+\int_{\mathbb{R}^N}\int_{\mathbb{R}^N}\frac{U^{2_\mu^*-2}(y)\psi^2(y)\partial_s\left(U^{2_\mu^*-2}(x)\psi^2(x)\right)}{|x-y|^{\mu}}dxdy
\Bigg)
=o\left(\lambda_{\varepsilon}^{N-2}\varepsilon^{N-2}
+\frac{1}{\lambda_{\varepsilon}^{N-2}}\right)
=o\left(\varepsilon^{\frac{N-2}{2}}
\right).
\end{align*}
In conclusion,
\begin{align}\label{C'0}
C'_0=\frac{N(N-2)^2}{\mathcal{A}_{H,L}}\left((N-2)\omega_{N}B_N\frac{H(0,0)\varepsilon^{-\theta}}{\lambda_{\varepsilon}^{N-1}}
-\frac{M(\tau)}{(1+|\tau|^2)^{\frac{N-2}{2}}}
\varepsilon^{N-2}\lambda_{\varepsilon}^{N-3}\varepsilon^{-\theta}\left(1+o(1)\right)
\right).
\end{align}
We now estimate $C'_1$, that is
\begin{align}\label{C'1}
C'_1=&O\left(\int_{\mathbb{R}^N\setminus \Omega_\varepsilon}\int_{\Omega_\varepsilon}
\frac{PU^{2_\mu^*}(y)PU^{2_\mu^*-1}(x)\partial_sPU(x)}{|x-y|^{\mu}}dxdy
\right)
\nonumber\\
=&O\left(\frac{N(N-2)}{\mathcal{A}_{H,L}}\int_{\mathbb{R}^N\setminus \Omega_\varepsilon}PU^{2_\mu^*-1}(x)\partial_sPU(x)dx
\right)
=O\left(\lambda_{\varepsilon}^{N}\varepsilon^{N}
+\frac{1}{\lambda_{\varepsilon}^{N}}\right)
=O\left(\varepsilon^{\frac{N}{2}}\right).
\end{align}
In the same way, we also obtain
\begin{align}\label{C'2}
C'_2=O\left(\lambda_{\varepsilon}^{N}\varepsilon^{N}
+\frac{1}{\lambda_{\varepsilon}^{N}}\right)
=O\left(\varepsilon^{\frac{N}{2}}\right).
\end{align}
Therefore, it follows from (\ref{C'0}), (\ref{C'1}) and (\ref{C'2}) that
\begin{align}\label{C'}
C'=\frac{N(N-2)^2}{\mathcal{A}_{H,L}}\left((N-2)\omega_{N}B_N\frac{H(0,0)\varepsilon^{-\theta}}{\lambda_{\varepsilon}^{N-1}}
-\frac{M(\tau)}{(1+|\tau|^2)^{\frac{N-2}{2}}}
\varepsilon^{N-2}\lambda_{\varepsilon}^{N-3}\varepsilon^{-\theta}\left(1+o(1)\right)
\right).
\end{align}
Combining (\ref{A'}), (\ref{B'}) and (\ref{C'}), we have
\begin{align*}
&\partial_{\lambda}I_\varepsilon(PU)
\nonumber\\
=&\frac{N(N-2)^3\omega_{N}B_{N}}{2\mathcal{A}_{H,L}\lambda_{\varepsilon}^{N-1}}H(0,0)\varepsilon^{-\theta}(1+o(1))
-\frac{N(N-2)^2}{2\mathcal{A}_{H,L}}\frac{M(\tau)}{(1+|\tau|^2)^{\frac{N-2}{2}}}
\varepsilon^{N-2}\lambda_{\varepsilon}^{N-3}\varepsilon^{-\theta}(1+o(1))
\nonumber\\
&-\frac{N(N-2)^2}{\mathcal{A}_{H,L}}\left((N-2)\omega_{N}B_N\frac{H(0,0)\varepsilon^{-\theta}}{\lambda_{\varepsilon}^{N-1}}
-\frac{M(\tau)}{(1+|\tau|^2)^{\frac{N-2}{2}}}
\varepsilon^{N-2}\lambda_{\varepsilon}^{N-3}\varepsilon^{-\theta}\left(1+o(1)\right)
\right)
\nonumber\\
=&-\frac{N(N-2)^3\omega_{N}B_{N}}{2\mathcal{A}_{H,L}\lambda_{\varepsilon}^{N-1}}
H(0,0)\varepsilon^{-\theta}(1+o(1))
+\frac{N(N-2)^2M(\tau)}{2\mathcal{A}_{H,L}(1+|\tau|^2)^{\frac{N-2}{2}}}
\varepsilon^{N-2}\lambda_{\varepsilon}^{N-3}\varepsilon^{-\theta}(1+o(1))
\nonumber\\
=&\frac{N(N-2)}{2\mathcal{A}_{H,L}}\varepsilon^{\frac{N-2}{2}}\partial_{\lambda}\Psi({\tau},{\lambda})
+o\left(\varepsilon^{\frac{N-2}{2}}\right).
\end{align*}
Similarly, for $\partial_s=\partial_{\tau^i}$ with $i=1,\cdots,N$, we have
\begin{align*}
\partial_{\tau^i}I_\varepsilon(PU)
=\frac{N(N-2)}{2\mathcal{A}_{H,L}}\varepsilon^{\frac{N-2}{2}}\partial_{\tau^i}\Psi({\tau},{\lambda})
+o\left(\varepsilon^{\frac{N-2}{2}}\right).
\end{align*}
That concludes the proof.
\end{proof}

In terms of Lemmas \ref{IPU} and \ref{IVphic1}, we have
\begin{align*}
I_\varepsilon(PU_{\lambda_\varepsilon, \xi_\varepsilon}+\phi)
=&\left(1-\frac{1}{2_\mu^*}\right)\frac{N(N-2)}{2\mathcal{A}_{H,L}}A_N
+\frac{N(N-2)^2\omega_{N}B_N}{2\mathcal{A}_{H,L}}\frac{H(0,0)}{\lambda_\varepsilon^{N-2}}\left(1+o(1)\right)
\nonumber\\
&+\frac{N(N-2)}{2\mathcal{A}_{H,L}(1+|\tau|^2)^{\frac{N-2}{2}}}
M(\tau)\varepsilon^{N-2}\lambda_\varepsilon^{N-2}\left(1+o(1)\right).
\end{align*}

Now we are ready to complete the proof of Proposition \ref{criticalpoint} (ii).

\noindent {\it Proof of Proposition \ref{criticalpoint} (ii)}:
It follows from Lemmas \ref{IVphic1} and \ref{IPU}.
\qed

\section{Appendix}\label{appendix}

\begin{lemma}\label{puipuj}
It holds
\begin{align*}
&\langle PZ_{\lambda_\varepsilon, \xi_\varepsilon}^0, PZ_{\lambda_\varepsilon, \xi_\varepsilon}^0\rangle=c_0+O\left(\frac{1}{\lambda_\varepsilon^{2N-\mu}}\right),\\
&\langle PZ_{\lambda_\varepsilon, \xi_\varepsilon}^0, PZ_{\lambda_\varepsilon, \xi_\varepsilon}^i\rangle=O\left(\frac{1}{\lambda_\varepsilon^{2N-\mu}}\right),\ i\neq  0,\\
&\langle PZ_{\lambda_\varepsilon, \xi_\varepsilon}^i, PZ_{\lambda_\varepsilon, \xi_\varepsilon}^j\rangle=\begin{cases}
o\left(\lambda_\varepsilon^2\right), &\ i\neq j,\\
c_{i}\lambda_\varepsilon^2+o\left(\lambda_\varepsilon^2\right), &\ i=j,
\end{cases}
\end{align*}
for some positive constants $c^i$ for $i=0,1,\cdots,N$.
\end{lemma}
\begin{proof}
By the definition of $PZ^i_{\lambda_\varepsilon,\xi_\varepsilon}$ and let $PZ^i_{\lambda_\varepsilon,\xi_\varepsilon}:=PZ^i$, $U_{\lambda_\varepsilon,\xi_\varepsilon}:=U$, we have
\begin{align*}
\langle PZ^i, PZ^j\rangle=\int_{\Omega_\varepsilon}f'(U)Z^iPZ^j
=\int_{\Omega_\varepsilon}f'(U)Z^iZ^j+\int_{\Omega_\varepsilon}f'(U)Z^i(PZ^j-Z^j),
\end{align*}
where
\begin{align*}
\int_{\Omega_\varepsilon}f'(U)Z^i(PZ^j-Z^j)
=O\left(\|Z^i\|\|PZ^j-Z^j\|\right)
=\begin{cases}
o(1), &\ i=0, \\
o(\lambda_\varepsilon^2), &\ i\neq 0,
\end{cases}
\end{align*}
and
\begin{align*}
\int_{\Omega_\varepsilon}f'(U)Z^iZ^j
=&2_\mu^*\int_{\Omega_\varepsilon}\int_{\Omega_\varepsilon}\frac{U^{2_\mu^*-1}Z^i(y)U^{2_\mu^*-1}Z^j(x)
}{|x-y|^{\mu}}dxdy
\nonumber\\
&+(2_\mu^*-1)\int_{\Omega_\varepsilon}\int_{\Omega_\varepsilon}\frac{U^{2_\mu^*}(y)U^{2_\mu^*-2}Z^jZ^i(x)
}{|x-y|^{\mu}}dxdy.
\end{align*}
Since
\begin{align*}
&2_\mu^*\int_{\Omega_\varepsilon}\int_{\Omega_\varepsilon}\frac{U^{2_\mu^*-1}Z^i(y)U^{2_\mu^*-1}Z^j(x)
}{|x-y|^{\mu}}dxdy
\nonumber\\
=&2_\mu^*(N-2)^2\int_{\Omega_\varepsilon}\int_{\Omega_\varepsilon}\frac{\left(\frac{\lambda_\varepsilon}{1+\lambda_\varepsilon^2|y-\xi_\varepsilon|^2}
\right)^{\frac{N-2}{2}(2_\mu^*-1)}
\frac{\lambda_\varepsilon^{\frac{N+2}{2}}(y-\xi_\varepsilon)_i}{\left(1+\lambda_\varepsilon^2|y-\xi_\varepsilon|^2\right)^{\frac{N}{2}}}
\left(\frac{\lambda_\varepsilon}{1+\lambda_\varepsilon^2|x-\xi_\varepsilon|^2}
\right)^{\frac{N-2}{2}(2_\mu^*-1)}
\frac{\lambda_\varepsilon^{\frac{N+2}{2}}(x-\xi_\varepsilon)_j}{\left(1+\lambda_\varepsilon^2|x-\xi_\varepsilon|^2\right)^{\frac{N}{2}}}
}{|x-y|^{\mu}}dxdy
\nonumber\\
=&2_\mu^*(N-2)^2\lambda_\varepsilon^2\int_{\lambda_\varepsilon(\Omega_\varepsilon-\xi_\varepsilon)}
\int_{\lambda_\varepsilon(\Omega_\varepsilon-\xi_\varepsilon)}\frac{
\frac{s_i}{\left(1+s^2\right)^{\frac{N-2}{2}(2_\mu^*-1)+\frac{N}{2}}}
\frac{t_j}{\left(1+t^2\right)^{\frac{N-2}{2}(2_\mu^*-1)+\frac{N}{2}}}
}{|s-t|^{\mu}}dsdt
\nonumber\\
=&\begin{cases}
o\left(\lambda_\varepsilon^2\right), &\ i\neq j,\\
c_{1i}\lambda_\varepsilon^2+o\left(\lambda_\varepsilon^2\right), &\ i=j,
\end{cases}
\end{align*}
and similarly,
\begin{align*}
&(2_\mu^*-1)\int_{\Omega_\varepsilon}\int_{\Omega_\varepsilon}\frac{U^{2_\mu^*}(y)U^{2_\mu^*-2}Z^jZ^i(x)
}{|x-y|^{\mu}}dxdy
\nonumber\\
=&(2_\mu^*-1)(N-2)^2\lambda_\varepsilon^2\int_{\lambda_\varepsilon(\Omega_\varepsilon-\xi_\varepsilon)}
\int_{\lambda_\varepsilon(\Omega_\varepsilon-\xi_\varepsilon)}\frac{
\frac{1}{\left(1+s^2\right)^{\frac{N-2}{2}(2_\mu^*)}}
\frac{t_it_j}{\left(1+t^2\right)^{\frac{N-2}{2}(2_\mu^*-2)+N}}
}{|s-t|^{\mu}}dsdt
\nonumber\\
=&\begin{cases}
o\left(\lambda_\varepsilon^2\right), &\ i\neq j,\\
c_{2i}\lambda_\varepsilon^2+o\left(\lambda_\varepsilon^2\right), &\ i=j,
\end{cases}
\end{align*}
we obtain that
\begin{align*}
\langle PZ^i, PZ^j\rangle=&
\begin{cases}
o\left(\lambda_\varepsilon^2\right), &\ i\neq j,\\
(c_{1i}+c_{2i})\lambda_\varepsilon^2+o\left(\lambda_\varepsilon^2\right), &\ i=j,
\end{cases}
:=
\begin{cases}
o\left(\lambda_\varepsilon^2\right), &\ i\neq j,\\
c_{i}\lambda_\varepsilon^2+o\left(\lambda_\varepsilon^2\right), &\ i=j.
\end{cases}
\end{align*}
In the same way, we get
\begin{align*}
&2_\mu^*\int_{\Omega_\varepsilon}\int_{\Omega_\varepsilon}\frac{U^{2_\mu^*-1}Z^i(y)U^{2_\mu^*-1}Z^0(x)
}{|x-y|^{\mu}}dxdy
\nonumber\\
=&2_\mu^*\frac{(N-2)^2}{2}\int_{\Omega_\varepsilon}\int_{\Omega_\varepsilon}\frac{\left(\frac{\lambda_\varepsilon}
{1+\lambda_\varepsilon^2|y-\xi_\varepsilon|^2}
\right)^{\frac{N-2}{2}(2_\mu^*-1)}
\frac{\lambda_\varepsilon^{\frac{N+2}{2}}(y-\xi_\varepsilon)_i}{\left(1+\lambda_\varepsilon^2|y-\xi_\varepsilon|^2\right)^{\frac{N}{2}}}
\left(\frac{\lambda_\varepsilon}{1+\lambda_\varepsilon^2|x-\xi_\varepsilon|^2}
\right)^{\frac{N-2}{2}(2_\mu^*-1)}
\lambda_\varepsilon^{\frac{N-4}{2}}\frac{1-\lambda_\varepsilon^2|x-\xi_\varepsilon|^2}{\left(1+\lambda_\varepsilon^2|x-\xi_\varepsilon|^2\right)^{\frac{N}{2}}}
}{|x-y|^{\mu}}dxdy
\nonumber\\
=&2_\mu^*\frac{(N-2)^2}{2}\int_{\lambda_\varepsilon(\Omega_\varepsilon-\xi_\varepsilon)}
\int_{\lambda_\varepsilon(\Omega_\varepsilon-\xi_\varepsilon)}\frac{
\frac{s_i}{\left(1+s^2\right)^{\frac{N-2}{2}(2_\mu^*-1)+\frac{N}{2}}}
\frac{1-t^2}{\left(1+t^2\right)^{\frac{N-2}{2}(2_\mu^*-1)+\frac{N}{2}}}
}{|s-t|^{\mu}}dsdt
=O\left(\frac{1}{\lambda_\varepsilon^{2N-\mu}}\right).
\end{align*}
Similarly,
\begin{align*}
(2_\mu^*-1)\int_{\Omega_\varepsilon}\int_{\Omega_\varepsilon}\frac{U^{2_\mu^*}(y)U^{2_\mu^*-2}Z^0Z^i(x)
}{|x-y|^{\mu}}dxdy
=O\left(\frac{1}{\lambda_\varepsilon^{2N-\mu}}\right),
\end{align*}
\begin{align*}
&2_\mu^*\int_{\Omega_\varepsilon}\int_{\Omega_\varepsilon}\frac{U^{2_\mu^*-1}Z^0(y)U^{2_\mu^*-1}Z^0(x)
}{|x-y|^{\mu}}dxdy
\nonumber\\
=&2_\mu^*\int_{\Omega_\varepsilon}\int_{\Omega_\varepsilon}\frac{\left(\frac{\lambda_\varepsilon}
{1+\lambda_\varepsilon^2|y-\xi_\varepsilon|^2}
\right)^{\frac{N-2}{2}(2_\mu^*-1)}
\frac{N-2}{2}\lambda_\varepsilon^{\frac{N-4}{2}}\frac{1-\lambda_\varepsilon^2|y-\xi_\varepsilon|^2}
{\left(1+\lambda_\varepsilon^2|y-\xi_\varepsilon|^2\right)^{\frac{N}{2}}}
\left(\frac{\lambda_\varepsilon}{1+\lambda_\varepsilon^2|x-\xi_\varepsilon|^2}
\right)^{\frac{N-2}{2}(2_\mu^*-1)}
\frac{N-2}{2}\lambda_\varepsilon^{\frac{N-4}{2}}\frac{1-\lambda_\varepsilon^2|x-\xi_\varepsilon|^2}
{\left(1+\lambda_\varepsilon^2|x-\xi_\varepsilon|^2\right)^{\frac{N}{2}}}
}{|x-y|^{\mu}}dxdy
\nonumber\\
=&2_\mu^*\frac{(N-2)^2}{4}\int_{\lambda_\varepsilon(\Omega_\varepsilon-\xi_\varepsilon)}
\int_{\lambda_\varepsilon(\Omega_\varepsilon-\xi_\varepsilon)}\frac{
\frac{1-s^2}{\left(1+s^2\right)^{\frac{N-2}{2}(2_\mu^*-1)+\frac{N}{2}}}
\frac{1-t^2}{\left(1+t^2\right)^{\frac{N-2}{2}(2_\mu^*-1)+\frac{N}{2}}}
}{|s-t|^{\mu}}dsdt
=c_{10}+O\left(\frac{1}{\lambda_\varepsilon^{2N-\mu}}\right),
\end{align*}
and
\begin{align*}
(2_\mu^*-1)\int_{\Omega_\varepsilon}\int_{\Omega_\varepsilon}\frac{U^{2_\mu^*}(y)U^{2_\mu^*-2}Z^0Z^0(x)
}{|x-y|^{\mu}}dxdy
=c_{20}+O\left(\frac{1}{\lambda_\varepsilon^{2N-\mu}}\right).
\end{align*}
In conclusion, we can conclude the proof.
\end{proof}

\subsection*{Acknowledgments}

Xiaomeng Huang was supported by the China Scholarship Council (No.202306990081). Marco Ghimenti is partially supported by the MIUR Excellence Department Project awarded to the Department of Mathematics, University of Pisa, CUP I57G22000700001 and by the
PRIN 2022 project 2022R537CS
\uline{$NO^3$-Nodal Optimization, NOnlinear elliptic equations,
NOnlocal geometric problems, with a focus on regularity}, founded by the European Union -Next Generation EU. Angela Pistoia is partially supported by the MUR-PRIN-20227HX33Z
``Pattern formation in nonlinear phenomena". The research of Marco Ghimenti and Angela Pistoia is partially supported by the GNAMPA project 2024: ``Problemi di doppia curvatura
su variet\`{a} a bordo e legami con le EDP di tipo ellittico"

\end{document}